\documentclass[11pt]{article}

\usepackage{amsmath, amsfonts,amssymb, amsthm, euscript,makeidx,color,mathrsfs,latexsym}

\newtheorem{assumption}{Assumption}

\newcommand{\la}{\langle}
\newcommand{\ra}{\rangle}
\newcommand{\hP}{\hat\dbP}

\newcommand{\ba}{\begin{array}}
\newcommand{\ea}{\end{array}}
\newcommand{\be}{\begin{equation}}
\newcommand{\ee}{\end{equation}}
\newcommand{\bee}{\begin{equation*}}
\newcommand{\eee}{\end{equation*}}
\newcommand{\bea}{\begin{eqnarray}}
\newcommand{\eea}{\end{eqnarray}}
\newcommand{\beaa}{\begin{eqnarray*}}
\newcommand{\eeaa}{\end{eqnarray*}}

\def\neg{\negthinspace}

\def\dbE{\mathbb{E}}
\def\dbF{\mathbb{F}}

\def\dbL{\mathbb{L}}

\def\dbP{\mathbb{P}}
\def\dbR{\mathbb{R}}
\def\dbS{\mathbb{S}}

\def\a{\alpha}
\def\b{\beta}
\def\g{\gamma}
\def\d{\delta}
\def\e{\varepsilon}

\def\si{\sigma}
\def\t{\tau}
\def\f{\varphi}
\def\th{\theta}
\def\o{\omega}

\def\G{\Gamma}
\def\D{\Delta}
\def\Th{\Theta}
\def\L{\Lambda}

\def\O{\Omega}
\def\cB{{\cal B}}

\def\cF{{\cal F}}
\def\cG{{\cal G}}

\def\cM{{\cal M}}
\def\cN{{\cal N}}

\def\cY{{\cal Y}}
\def\cZ{{\cal Z}}
\def\hE{\mathbb{E}}
\def\hF{\mathbb{F}}

\def\hL{\mathbb{L}}

\def\hN{\mathbb{N}}

\def\hP{\mathbb{P}}

\def\hR{\mathbb{R}}
\def\hS{\mathbb{S}}

\def\hX{\mathbb{X}}

\def\sA{\mathscr{A}}
\def\sB{\mathscr{B}}

\def\sD{\mathscr{D}}

\def\sM{\mathscr{M}}
\def\sN{\mathscr{N}}

\def\sU{\mathscr{U}}

\def\sY{\mathscr{Y}}
\def\sZ{\mathscr{Z}}
\def\no{\noindent}

\def\ms{\medskip}
\def\bs{\bigskip}
\def\q{\quad}
\def\qq{\qquad}

\def\pa{\partial}
\def\cd{\cdot}
\def\cds{\cdots}

\def\tr{\hbox{\rm tr}}
\def\ol{\overline}

\newcommand{\basa}{\begin{assumption}}
\newcommand{\easa}{\end{assumption}}

\newcommand{\bas}{\begin{assum}}
\newcommand{\eas}{\end{assum}}

\def\essinf{\mathop{\rm essinf}}
\def\esssup{\mathop{\rm esssup}}

\def\pa{\partial}

 \def\cd{\cdot}
\def\cds{\cdots}

\def\as{\hbox{\rm -a.s.{ }}}

\def\tr{\hbox{\rm tr$\,$}}

\def\dis{\displaystyle}

\def\1{{\bf 1}}

\def\:{\!:\!}
\def\reff#1{{\rm(\ref{#1})}}
\def \proof{{\noindent \bf Proof\quad}}

at 9pt

\begin{document}

\newtheorem{thm}{Theorem}[section]
\newtheorem{lem}[thm]{Lemma}
\newtheorem{cor}[thm]{Corollary}
\newtheorem{prop}[thm]{Proposition}
\newtheorem{rem}[thm]{Remark}
\newtheorem{eg}[thm]{Example}
\newtheorem{defn}[thm]{Definition}
\newtheorem{assum}[thm]{Assumption}

\renewcommand {\theequation}{\arabic{section}.\arabic{equation}}
\def\thesection{\arabic{section}}

\title{\bf  Dynamic Approaches for Some Time Inconsistent Problems}

\author{
Chandrasekhar Karnam\thanks{ \no Department of
Mathematics, University of Southern California, Los Angeles, 90089;
email: chandrasekhar.karnam@gmail.com},  ~ Jin Ma\thanks{ \noindent Department of
Mathematics, University of Southern California, Los Angeles, 90089;
email: jinma@usc.edu. This author is supported in part
by US NSF grant \#1106853. } ~ and ~{Jianfeng Zhang}\thanks{\noindent
Department of Mathematics, University of Southern California, Los
Angeles, CA 90089. E-mail: jianfenz@usc.edu. This author is
supported in part by NSF grant \#1413717. 
}}

\date{\today}
\maketitle

\begin{abstract}
In this paper we investigate possible approaches to study general time-inconsistent optimization problems without assuming
the existence of optimal strategy. This leads
immediately to the need to refine the concept of time-consistency as well as any method that is based on  Pontryagin's Maximum Principle. The fundamental obstacle is  the dilemma of having to invoke the {\it Dynamic Programming Principle} (DPP) in a
time-inconsistent setting, which is contradictory in nature. The main contribution of this work is the introduction of the idea 
of the ``dynamic utility" under which the original time inconsistent problem (under the fixed utility) becomes a time consistent one. As a benchmark model, we shall consider a stochastic controlled problem with multidimensional backward SDE
dynamics, which covers many existing time-inconsistent problems in the literature as special cases; and we argue that the time inconsistency is essentially equivalent to the lack of {\it comparison principle}.
 We shall propose three approaches aiming  at reviving the DPP in this setting: the duality approach, the dynamic utility approach, and the master equation approach. 
 Unlike the game approach in many existing works in continuous time models,  all our approaches produce the same value as the  original static problem. 
\end{abstract}

\vfill \bs

\no


\bs

\no{\it 2000 AMS Mathematics subject classification:}  	49L20, 60H10,   91C99, 91G80, 35R15

\eject

\section{Introduction}
\label{sect-Introduction}
\setcounter{equation}{0}

In this paper we propose some possible approaches to tackle the general {\it time-inconsistent} optimization
problems in continuous time setting. These approaches are different from all the existing ones in the literature, and are based on our
new understanding of the time inconsistency.
We note that the time inconsistency appears naturally and frequently in economics 
and finance,
 see e.g. Kydland-Prescott \cite{KP} and  Kahneman-Tversky \cite{KT1, KT2}. We refer to the frequently cited survey Strotz \cite{Strotz} for the fundamentals of this problem, and Zhou \cite{Zhou} for some recent development on continuous time models. 
We should point out that it was 
 \cite{Zhou} that brought  the time inconsistency issue to our attention.

\ms
\no{\bf I. Time inconsistency.} We begin by briefly describing the time-inconsistency in an optimization problem that has been understood so far. Consider an optimization problem over a time interval $[0, T]$:
\bea
\label{intro-V0}
V_0 := \sup_{u\in \sU_{[0,T]}} J(u).
\eea
where $\sU_{[0,T]}$ is an appropriate set of admissible controls $u$ defined on $[0,T]$, and $J(u)$ is a certain 
utility functional  associated to  $u$. Clearly, the problem \reff{intro-V0} is static.   Its dynamic counterpart is the following optimization problem over $[t, T]$, for any $t\in [0, T]$:
\bea
\label{intro-Vt}
V_t := \esssup_{u\in \sU_{[t, T]}} J_t(u).
\eea
Here $\sU_{[t, T]}$ is the corresponding set of admissible controls on $[t, T]$   and utility functional $J_t$ usually involves some conditional expectation, and thus could be random. 

 An admissible control $u^*\in\sU_{[0,T]}$ is called ``optimal"  for the  problem \reff{intro-V0} if $J(u^*)=V_0$. Defining optimal control $u^{t,*}$ for the problem \reff{intro-Vt} similarly and assuming their existence, we say the problem \reff{intro-Vt}  is {\it time-consistent} if, for any $t\in [0, T]$,  it holds that
\bea
\label{consistent}
u^{t, *}_s = u^*_s, \q t\le s\le T.
\eea
The relation (\ref{consistent}) amounts to saying that   {\it a (temporally) global optimum must be a 
local one}. 
The optimization problem \reff{intro-Vt} is called {\it time-inconsistent} if \reff{consistent} fails to hold. Intuitively, time inconsistency 
means an optimal strategy today may not be optimal tomorrow.

Since the early work \cite{Strotz}, there have been typically two approaches for treating the time inconsistent problems, both focusing on the optimal control: (i) the strategy of {\it precommitment}, and (ii) the strategy of {\it consistent planning}.   The former is to solve the static optimization problem \reff{intro-V0}, and then simply insist on using $u^*$ (assuming it exists) throughout $[0, T]$, despite the fact that it may not be optimal anymore when $t>0$. The latter one has
developed into the popular ``game approach" in the literature, in which the player plays with infinitely many future selves.  To illustrate the idea, let us consider the discrete time setting: $0=t_0<\cds<t_n=T$. The ``consistent planning" amounts to saying that at any $t_i$, the player tries to find optimal strategy $u$ on $[t_i, t_{i+1})$ by assuming the future selves have already found the optimal strategies and will actually use them on $[t_{i+1}, T]=[t_{i+1}, t_{i+2})\cup\cds\cup[t_{n-1},T]$.  We note that an equilibrium in such a game approach should
be similar to that of a principal agent problem, that is, in the sense of a  sequential optimization problem, rather than a Nash equilibrium. 

The game approach makes sense in many applications, but is very challenging in continuous time setting (being a game with uncountably many players!). There have been some successful applications of this approach in continuous time models, see, 
e.g., 
 Bjork \& Murgoci \cite{BM},  Ekeland \& Lazrak \cite{EL}, Hu, Jin \& Zhou \cite{HJZ}, and Yong  \cite{Yong}, to mention a few.
It is worth noting that since under the game framework the problem is time consistent, which enables one to apply the standard tools such as dynamic programming and HJB equations. However, typically the value of the game problem at 
$t=0$ is different from the original value $V_0$ in \reff{intro-V0} (unless the problem is time consistent),  thus the solution of the game approach, even if it exists, does not really solve the 
problem \reff{intro-V0}.

In this paper we will be  focusing on 
the value $V_0$ of the original static problem \reff{intro-V0}. 
We would like to emphasize that the problem \reff{intro-V0}, or its ``precommitment" nature, actually makes more sense in some applications. For example, in the so-called 
{\it principal-agent problem} (see \S\ref{sect-PA}  below), practically the principal cannot change the contract once it commenced (at least not as frequently as the game approach requires), 
therefore one is obliged to follow the contract designed at $t=0$ for the whole contractual period. In fact, 
problem \reff{intro-V0} 
is a mathematically interesting problem in its own right. 

Another main feature of this paper is that, unlike most of the works in the ``time inconsistency" literature to date, we shall remove the presumption of the existence of optimal strategy. In fact, as is well known in stochastic control literature, it is not unusual that the optimal control fail to exist. 
It is somewhat surprising that without the
optimal control (or equilibrium in game approach), it is even not clear how to define the notion of time 
consistency/inconsistency(!) in the current literature, much less to 
say anything about the value $V_0$, which on the other hand is always well defined, regardless the existence of optimal control. Our main task is thus to find the new (time consistent) methods to solve the original value $V_0$, and to revive the dynamical programming method in a novel context. 

\ms
\no{\bf II. Our main observation.} It is well-understood that there are typically two approaches to solve the optimization problem \reff{intro-V0}: the {\it Dynamic Programming Principle} (DPP for short) and the {\it Stochastic Maximum Principle} (SMP for short). The former
relies fundamentally on the time consistency; whereas the latter
requires the existence of optimal control.  We then immediately find ourselves facing the dilemma: on the one hand  the SMP, as 
a {\it necessary condition},  is no longer relevant without an optimal control; 
but on the other hand,  DPP does not make sense either due to the lack of time-consistency.

To ``revive" the DPP for the static problem \reff{intro-V0}, our first plan is based on the following simple but crucial observation: the problem \reff{intro-Vt}  is time inconsistent partly due to the fact that, 
modulus some conditional expectation, the utility  $J_t$ in \reff{intro-Vt} is essentially the same as the utility $J$ in \reff{intro-V0}, which could be
in conflict with the nature of the problem and causing the time inconsistency. Therefore, if we allow $J_t$ to vary more freely 
with the time $t$, denoting it by $J(t, u)$, then it is hopeful that the new dynamic optimization problem
\bea
\label{intro-tildeVt}
\tilde V_t := \sup_{u\in \sU_{[t, T]}} J(t, u)
\eea
could become time consistent with the right choice of $J(t,\cd)$. In particular,  if we require that $J(0, u) = J(u)$, then $\tilde V_0 = V_0$ and  we are indeed solving the original problem \reff{intro-V0}. In fact, as we will see  in the next section, when the optimal control $u^*$ exists, one can easily construct such $J(t,\cd)$ by utilizing the optimal $u^*$.  The real challenge is, of course, to find a desired $J(t,\cd)$ without using $u^*$ or  in the situation where $u^*$ does not exist. 

We remark that,  given the initial value $J(0, u) = J(u)$,  the dynamic $J(\cd,\cd)$ will be sought {\it forwardly} (in time), and thus it is in spirit similar to the notion of {\it forward utility} proposed in
 \cite{MZ1, MZ2}.
 However, it should be emphasized that the dynamic utility $U(t,\cd)$ in \cite{MZ1, MZ2} is applied on an optimization problem over time period $[0, t]$, while our dynamic $J(t, \cd)$ is over time period $[t, T]$. Namely, there is a fundamental difference between the two notions.

We remark that  similar ideas of such ``dynamic utilities" have also appeared in the literature  in various different contexts, see for example,
Bouchard,  Elie \& Touzi \cite{BET}, Cui, Li, Wang \& Zhu \cite{CLWZ}, and Feinstein \& Rudloff \cite{FR2}. 

\ms
\no {\bf III. Our proposed approaches.} Our second main observation in this paper is that many time inconsistent problems in the literature can be transformed into control problems on multidimensional (possibly infinite dimensional) forward-backward SDEs
(see \S2 for details). Therefore in what follows we shall focus on the following benchmark optimization problem for controlled multidimensional backward SDEs:
\bea
\label{BSDE-V0}
 V_0 := \sup_{u\in \sU_{[0,T]}} \f(Y^u_0) ~\mbox{where}~ Y^u_t = \xi + \int_t^T f(s, Y^u_s, Z^u_s, u_s) ds - \int_t^T Z^u_s dB_s, ~ t\in[0,T].
\eea
 We note that in (\ref{BSDE-V0}) we have made two simplifications in order to focus more on the main issue of time inconsistency:  the controlled dynamics is only a backward SDE
 and the dimension is finite. All the results in this paper can be extended to the controlled forward-backward SDE case, but with
 heavier presentations. We prefer not to seek such generality in this paper. 
 The infinite dimensional case, however,  is more challenging, and we shall leave it to future study.

We start with a ``duality approach" by first noticing that
\bea
\label{intro-cD0}
V_0 = \sup_{y\in \sD_0} \f(y) \q\mbox{where}\q \sD_0 := \{Y^u_0: u\in \sU_{[0, T]}\}.
\eea
We shall argue that, in the Markovian case, the ``reachable set" $\sD_0$ can be written as 
\bea
\label{intro-N0}
\ol \sD_0 =  \sN(0,0):= \{y: W(0,0, y) = 0\},
\eea
where $\ol \sD_0$ is the closure of $\sD_0$, $W(t,x,y)$ is the unique viscosity solution to certain standard HJB equation, and $\sN(0,0)$ is the the so-called ``nodal set" of $W$. Assuming $\f$ is continuous, we can first solve the HJB equation for $W$, then compute its nodal set $\sN(0,0)$, and finally solve a simple {\it finite dimensional} optimization problem: 
\bea
\label{intro-V0W}
V_0 = \sup_{y\in \sN(0,0)} \f(y).
\eea
We note that the idea of nodal set was used in Ma \& Yong \cite{MY} for solving a forward-backward SDE (without control $u$),
and we call this a ``duality approach". We shall further argue that 
the duality  holds  in non-Markovian case as well, by utilizing the viscosity theory of path dependent PDEs developed  by Ekren, Keller, Touzi \& Zhang \cite{EKTZ} and
Ekren, Touzi \& Zhang \cite{ETZ1, ETZ2}.

While the duality approach is quite generally applicable under mild conditions,  it solves only the static problem $V_0$. In particular, it does not provide a time consistent dynamic value $\tilde V_t$.  Our next step is to extend the set $\sD_0$ and the duality  \reff{intro-N0} to a dynamic version:
\bea
\label{intro-cDt}
 \sD_t := \{Y^u_t: u\in \sU_{[t, T]}\},\q \ol {\sD_t} = \sN(t, B_t):= \{y: W(t, B_t, y) = 0\}.
\eea
We shall argue that the family $\{\sD_t\}_{0\le t\le T}$ satisfies a  {\it geometric DPP}, in the spirit of 
Soner \& Touzi \cite{ST},  and closely related to the set valued analysis (see e.g. 
Aubin \&  Frankowska \cite{AF} and Feinstein \& Rudloff \cite{FR1}). 
However, we note that the following natural dynamic value 
\bea
\label{BSDE-Vt}
V_t := \esssup_{u\in \sU_{[t,T]} }\f(Y^u_t) = \esssup_{y\in \sD_t} \f(y) = \esssup_{y\in \sN(t, B_t)} \f(y)
\eea
is typically time inconsistent. Here $\esssup_{y\in \sD_t} \f(y)$ means $\esssup_{y\in \dbR^d} [\f(y) \1_{\sD_t}(y)]$, the same for other similar notations.  
The goal of our second approach
is to find a {\it dynamic utility function} $\Phi(t, y)$ (possibly random) satisfying $\Phi(0,\cd) = \f$ and that
\bea
\label{BSDE-tildeVt}
\tilde V_t :=  \esssup_{u\in \sU_{[t,T]} }\Phi(t, Y^u_t)= \esssup_{y\in \sD_t} \Phi(t,y) = \esssup_{y\in \sN(t, B_t)} \Phi(t,y)
\eea
is time consistent. We shall name this the ``dynamic utility approach" for simplicity. An important observation coming out 
from the study of this approach is that the time inconsistency of \reff{BSDE-Vt} is essentially equivalent to the lack of comparison principle for the multidimensional  BSDE, a well-known fact in BSDE theory. Thus our task becomes to 
find some dynamic utility function $\Phi(t,\cd)$ which satisfies a certain comparison principle. In this paper we
succeed in finding a desired $\Phi$ in a linear case, and we shall leave the general nonlinear case, which seems to be 
quite challenging, 
to  future research.

Our last approach 
borrows the idea from the mean field game literature (see e.g. Cardaliaguet, Delarue,  Lasry \& Lions \cite{CDLL}), which we now describe.
First note that  the value $V_0$ in \reff{intro-cD0} is clearly 
a function of  terminal condition $\xi$. Thus, for any $t\in [0, T]$ and random variable $\eta \in \dbL^2(\cF_t)$, we define
\bea
\label{intro-Psi}
\Psi(t, \eta) := \sup_{u\in \sU_{[0, t]}} \f (\sY^u_0(t,\eta)),
\eea
where $\sY^u(t,\eta)$ is the solution to BSDE  \reff{BSDE-V0} on $[0, t]$, satisfying $\sY^u_t(t,\eta) = \eta$. Clearly,  $\Psi(0, y) = \f(y)$ and $V_0 = \Psi(T,\xi)$, thus both functions $\Phi$ in \reff{BSDE-tildeVt} and  $\Psi$ in \reff{intro-Psi} are temporally ``dynamic" in nature, with the same initial value $\f$. The main difference, however, is that in \reff{BSDE-tildeVt} the control is over $[t, T]$, whereas in \reff{intro-Psi} the control is over $[0, t]$. One should also note that, unlike in mean field theory where the functions often depend only on the laws of the random variables, the function $\Psi$ in (\ref{intro-Psi}) depends indeed on the random variable $\eta$, or more precisely on the joint law of $(\eta, B)$. 

A very pleasant surprise of the (forward) value function $\Psi$ is that it satisfies the following form of DPP almost automatically,
and can thus be viewed as  time consistent: 
\bea
\label{intro-PsiDPP}
\Psi(t_2, \eta) := \sup_{u\in \sU_{[t_1, t_2]}} \Psi (t_1, \sY^u_{t_1}(t_2,\eta)),~ \eta \in \dbL^2(\cF_{t_2}),\q \mbox{for any}~ 0\le t_1 < t_2 \le T.
\eea
We shall emphasize that, unlike the usual DPP in stochastic control literature, (\ref{intro-PsiDPP}) is {\it forward} (in time), i.e., 
$t_1<t_2$(!).
This is due to the fact that we are optimizing a {\it backward} controlled  problem.   To the best of our knowledge, such type of forward DPP is new.

Having obtained the DPP \reff{intro-PsiDPP}, we believe that certain HJB type of differential equation (for $\Psi$)
should naturally come into
the picture, which we shall name as the {\it master equation}, due to  the nature of  the function $\Psi$.  
We expect  two features for this master equation: first, it should be a first order partial differential equation in a certain sense, due to the forward nature of the DPP; second, 
it should involve certain path derivatives of  $\eta$ 
in the sense of Dupire \cite{Dupire}, due to the progressive measurability of  $\Psi$ and the requirement 
 $\eta$ being $\cF_t$-measurable. 
We shall argue that when the function $\Psi$ defined by \reff{intro-Psi} is smooth (to be specified in the paper), it will be the unique (classical) solution to our master equation. The main
difficulty of this approach, however, is when $\Psi$ does not have the desired smoothness. It then becomes a very interesting, albeit challenging, problem to propose appropriate notion of weaker solution to the master equation.
We shall leave this to future research.

The rest of the paper is organized as follows. In \S\ref{sect-inconsistent} we present several examples of time inconsistent problems. 
 In \S\ref{sect-BSDE} we introduce our model and explain the role of comparison principle in time consistency issue. In \S\ref{sect-duality}-\ref{sect-master} we propose the three approaches, respectively. 

\section{Preliminaries and Examples}
\label{sect-inconsistent}
\setcounter{equation}{0}

Throughout this paper we shall use the following canonical setup. Let $T>0$ be a fixed time horizon, $\O:= \{\o\in C([0, T], \dbR^d): \o_0=0\}$ the canonical space, $\cF=\sB(\O)$, the Borel $\si$-filed of $\O$, and $\dbP_0$ the Wiener measure.
Further, we let $B_t(\o)=\o_t$, $\o\in\O$ be the canonical process and $\dbF := \dbF^B$ the natural filtration generated by $B$, augmented by $\hP_0$. Then $B$ is an $\hF$-Brownian motion under $\hP_0$.  We also denote $\dbE := \dbE^{\dbP_0}$ for simplicity, when the contact is clear, and $\L := [0, T]\times \O$.

For a generic Euclidean space $\hX$, we denote its inner product by $(x, y)=x \cd y=x^\top y$, its
norm by $|x|=(x,x)^{1/2}$, and its Borel $\si$-field by $\sB(\hX)$. If $\hX=\hR^{d_1\times d_2}$, we denote $A:B=\tr (A^\top B)$, 
for $A, B\in\hX$. Also, let $\cG\subseteq \cF$ be any sub-$\si$-field and
$[s,t]\subseteq [0,T]$,
we denote 

$\bullet$ $\dbL^2({\cG};\hX)$ to be all $\hX$-valued, $\cG$-measurable random variable $\xi$ such that 
$\|\xi\|^2_2:=  \dbE[|\xi|^2]<\infty$. The inner product in $\dbL^2(\cG;\hX)$ is denoted by $(\xi,\eta)_{2}:= \hE[(\xi, \eta)]$, 
$ \xi$, $\eta\in \dbL^2({\cal G}; \hX)$. 
 
$\bullet$ $\dbL^2_\dbF([s, t];\hX)$ to be all $\hX$-valued, $\dbF$-adapted process $\eta$  on $[s,t]$, such that 
$$\|\eta\|_{2,s,t}:=\dbE\big[\int_s^t|\eta_t|^2dt\big]^{1/2}<\infty;$$
In particular, if  $\hX=\hR$, we shall omit $\hX$ in the above notations for simplicity.

 
In what follows we present several examples of  time inconsistent optimization problems. In each of these examples we shall 
see the BSDE formulation of the original problem and the possibility of finding the dynamic utility. For simplicity, in this section we assume $d=1$.

\subsection{A mean-variance optimization problem}
\label{sect-mv}

Consider a simple controlled stochastic dynamics
\bea
\label{X0}
X^u_s= x_0 + \int_0^s u_r dr + \int_0^s u_r dB_r,  \q s\in[0,T], \q u\in \sU_{[0,T]} := L^2_\dbF([0,T]).
\eea
 Let  $c>0$ be a constant, and consider the optimization problem:
\bea
\label{mv-V0}
V_0 := \sup_{u\in\sU_{[0,T]}} \Big\{\dbE[X^u_T] - {1\over 2c} \mbox{Var}(X^u_T)\Big\}. 
\eea
Following the arguments in \cite{HJZ}, one shows that the above optimization problem has an optimal feedback control:
$
u^*(s,x) = x_0 - x+ c e^T$, $0\le s\le T$. In other words, the optimal control  is:
$u^*_s = u^*(s, X^*)=x_0 - X^*_s + c e^T$, $s\in[0,T]$, where $X^*$ is the corresponding optimal dynamics satisfying
\beaa
 X^*_s = x_0 + \int_0^s [x_0-X^*_r + c e^T] dr + \int_0^s [x_0-X^*_r + c e^T] dB_r, \q s\in[0,T].
\eeaa

Now let $0<t<T$ be given, and we follow the control $u^*$ on $[0, t]$ so that $X^*_t$ is well-defined.  Consider the optimization problem on $[t, T]$, 
starting from $X^*_t$:
\bea
\label{X1}
X^{t,u}_s = X^*_t + \int_t^s u_r dr + \int_t^s u_r dB_r, \q s\in [t, T];
\eea
and define, similar to (\ref{mv-V0}), the value of the optimization problem at time $t$:
\bea
\label{mv-Vt}
V_t := \esssup_{u\in \sU_{[t,T]}}\Big\{\dbE_t[X^{t,u}_T] - {1\over 2c} \mbox{Var}_t(X^{t,u}_T)\Big\},
\eea
where $\hE_t[\cd]=\hE[\cd|\cF_t]$, and Var$_t$ is the conditional variance under $\hE_t$. Again, as before we should have optimal control 
on $[t, T]$: $u^{t, *}(s,x) = X^*_t - x + c e^{T-t}$, $s\in[t,T]$. It is clear that $u^{t,*}(s,x) \neq u^*(s,x)$. Consequently,  $u^{t,*}_s := u^{t,*}(s,X^{t,*}_s) \neq u^*_s$, where
\beaa
X^{t,*}_s = X^*_t + \int_t^s [X^*_t - X^{t,*}_r + c e^{T-t}] dr + \int_t^s  [X^*_t - X^{t,*}_r + c e^{T-t}] dB_r, \q s\in [t, T];
\eeaa
  Thus the problem \reff{X1}-\reff{mv-Vt} is time inconsistent.

However, we should note that we can change the cost functional in \reff{mv-Vt} slightly so that it becomes time consistent. In fact, let $c_t>0$ be a random process and consider 
\bea
\label{mv-tildeVt}
\tilde V_t := \esssup_{u\in\sU_{[t,T]}}\Big\{\dbE_t[X^{t,u}_T] - {1\over 2 c_t} \mbox{Var}_t(X^{t,u}_T)\Big\}.
\eea
A similar argument would lead us to the optimal feedback control: 
$\tilde u^{t, *}(s,x) = X^*_t - x + c_t e^{T-t}$.
If we set 
\bea
\label{mv-ct}
 c_t := ce^t - e^{t-T}[X^*_t - x_0],\q t\in[0,T],
 \eea
  then $\tilde u^{t, *}(s,x)  = x_0 - x + c e^T = u^*(s,x)$. 
Namely the  problem \reff{X1} \&  \reff{mv-tildeVt} is time consistent. 
\begin{rem}
\label{remark2.1}{\rm
(i) Since $c_0=c$, we have $\tilde V_0=V_0$. To wit,  $\{\tilde V_t\}_{0\le t\le T}$ is a time consistent dynamic system with
 initial value $V_0$, as desired. 

(ii) We note that in the portfolio selection problems,  the constant $c$ in (\ref{mv-V0}) usually stands for the risk aversion parameter of the investor. In practice, it is reasonable that this risk aversion parameter  may evolve as time changes. 
A time inconsistent problem where the constant $c$ depends on  state process $X$ was studied in \cite{BMZ}. Our example shows that if $c_t$ is chosen correctly, then the problem could become time consistent. 

(iii)  A discrete case in the same spirit of this example was studied in \cite{CLWZ}. 
\qed
}
\end{rem}

It is worth noting that the parameter $c_t$ in \reff{mv-ct} is constructed via the optimal control $u^*$ (and so are the examples in \S\ref{sect-1d}, \ref{sect-PA}), which is  undesirable given our goal  of  tackling the time inconsistency without using optimal strategy. Such a slight drawback notwithstanding, an important observation from this example is that 
the problem \reff{X0}-\reff{mv-V0} can be converted to an 
optimal control problem for a 2-dimensional Backward SDE:
\bea
\label{mv-BSDE}
\left.\ba{lll}
\dis V_0 := \sup_{u\in \sU} \f(Y^{1,u}_0, Y^{2, u}_0),\q\mbox{where}\q \f(y_1, y_2) :=   y_1  + {1\over 2c} |y_1|^2 -{1\over 2c} y_2,\\
\dis Y^{1,u}_t= X^u_T -\int_t^T Z^{1,u}_sdB_s,\q Y^{2,u}_t= |X^u_T|^2 -\int_t^T Z^{2,u}_sdB_s, \q t\in[0,T].\\
%
\ea\right.
\eea
As we pointed out  in Introduction and will articulate more in next section, one of the main reasons for the time inconsistency 
is the lack of {\it comparison principle} for the underlying  dynamics, which is particularly the case for  \reff{mv-BSDE}. 

\subsection{A one dimensional example}
\label{sect-1d}

Besides the comparison principle as mentioned in the end of the previous subsection, another reason for time inconsistency is that the $\f$ in \reff{mv-BSDE} is not monotone. In what follows we present a one dimensional example where the comparison principle holds true. 

Let $\sU:=L^2_\dbF([0,T]; [-1,1])$.
Consider a simple one-dimensional BSDE:
\bea
\label{Y0}
Y^u_s= B_T + \int_s^T u_r dr - \int_s^T Z^u_r dB_r, \q s\in[0,T], \q u\in\sU,
\eea
and, let $\f(y) := -|c+y|$, $y\in\hR$, for some constant $c\in \dbR$. We define the optimal value by 
\bea
\label{1d-V0}
V_0 := \sup_{u\in \sU} \f(Y^u_0)= \sup_{u\in \sU} \f(\hE[Y^u_0])
= -\inf_{u\in \sU} \Big| c+ \int_0^T \dbE[u_s]ds \Big|.
\eea
Then one can easily check that
$u^*\in \sU$ is an optimal control if and only if:
\beaa
u^*_s \equiv -1, ~\mbox{if}~c\ge T;\q u^*_s\equiv1,~\mbox{if}~c\le -T;\q  \mbox{and}~\int_0^T \dbE[u_s]ds = -c, ~\mbox{if}~ |c|< T.
\eeaa

Now assume $c= T$. Let $0<t<T$ and 
consider the optimization problem over $[t, T]$:
\bea
\label{1d-Vt}
\dis V_t := \esssup_{u\in \sU} \f(Y^{u}_t) = -  \essinf_{u\in \sU} \Big|T+B_t + \int_t^T \dbE_t[u_s]ds\Big|, 
\eea
where $\hE_t[\cd]=\hE[\cd|\cF_t]$. Since $c=T$, if the problem were time-consistent we would then expect that the optimal control is $u^*_s=-1$, from the previous argument. However, 
we note that on the set $\{B_t \le t-2T\}$, one has 
$$0\ge T+B_t  +(T-t)\ge T+B_t+\int_t^T \hE_t[u_s]ds, \qq \mbox{for all $u\in\sU$},
$$ 
thus the optimal control for $V_t$ should be $u^{t, *}_s= 1$ on the set $\{B_t\le t-2T\}$, instead of $u^*_s=-1$, a contradiction. Namely the problem (\ref{1d-V0}) is time-inconsistent. 

Similar to the example in the previous subsection, if we allow the constant $c$ in (\ref{1d-V0}) to be time varying and even
random, then the problem could become time consistent. Indeed, if we choose $c_t := T-t- B_t$, and consider
\bea
\label{1d-tildeVt}
\tilde V_t := \esssup_{u\in \sU} \Phi(t, Y^{u}_t), \q\mbox{where}\q \Phi(t, y) := -|c_t + y|.
 \eea
Then it is readily seen that
\beaa
\tilde V_t = - \essinf_{u\in \sU} \Big|(T-t-B_t) + B_t + \int_t^T \dbE_t[u_s]ds\Big| = -  \essinf_{u\in \sU} \Big|T-t+ \int_t^T \dbE_t[u_s]ds\Big|,
\eeaa
and thus the optimal control is still $u^*=-1$.

\subsection{A principal-agent problem}
\label{sect-PA}

In this example we consider a special case of the Holmstrom-Milgrom model in the {\it Pringcipal-agent Problem} (cf. \cite{CZ}).
In this problem the principal is to find the optimal contract assuming the agent(s) will always perform optimally given any contract. 
The main feature of principal's contract is that it is pre-committed, that is, it cannot be changed (at least not frequently) during a contractually designed duration.

To be more precise, let $\g_A>0, \g_P>0$, $R<0$ be constants, and consider two exponential utility functions:
\beaa
U_A(x) := -\exp\{-\g_A x\},\q  U_P(x) := -\exp\{-\g_P x\}.
\eeaa
We denote the principal's control set by $\sU_{P}\subset \dbL^2(\cF_T)$, and the agent's control set by $\sU_A\subset \dbL^2_\dbF([0,T])$, satisfying certain technical conditions which for simplicity we will not specify.
 Given any contract $C_T\in \sU_P$ at $t=0$, we consider the agent's problem:
\bea
\label{PA-VA0}
 V^A_0(C_T) := \sup_{u\in\sU_A} \dbE^{\dbP^u}\Big[ U_A\big(C_T  - {1\over 2}\int_0^T |u_s|^2 ds\big)\Big], 
\eea
where $\hP^u$ is a new probability measure defined by $ {d\dbP^u\over d\dbP_0} := M^u_T := \exp\big\{\int_0^T u_s ds - {1\over 2}\int_0^T |u_s|^2 ds\big\}$.  We note that here the agent's control problem (\ref{PA-VA0}) is in a ``weak formulation", and 
$V_0^A(C_T)\le 0$ is well-defined. We shall consider those contracts that satisfy the following ``participation constraint"
\bea
\label{V0geR}
 V_0^A(C_T) \ge R, 
\eea
where $R<0$ is the ``market value" of an agent that a principle has to consider at $t=0$.

It can be shown  (cf. \cite[Chapter 6]{CZ}) that the agent's problem can be solved in terms of the following quadratic BSDE:
\beaa
Y^A_s = C_T - {\g_A-1\over 2}\int_s^T |Z^A_r|^2  dr - \int_s^T Z^A_r dB_r, \qq s\in [0,T].
\eeaa
In fact, by a simple comparison argument for BSDEs one shows that the  agent's optimal action is $u^*=u^*(C_T) = Z^A\in \sU_A$, with optimal value
$V^A_0 = U_A( Y^A_0)$.

Given the optimal $u^*=u^*(C_T)$ we now consider the principal's problem:
\bea
\label{PA-VP0}
V^P_0 := \sup_{C_T\in\sU_P} \dbE^{\dbP^{u^*}}[U_P(B_T-C_T)], 
\eea
subject to the participation constraint (\ref{V0geR}). The solution to the problem (\ref{PA-VP0})-(\ref{V0geR}) can be found explicitly
(cf. \cite[Chapter 6]{CZ}). Indeed, the optimal contract is:
\beaa
\dis C_T^* := - {1\over \g_A}\ln (-R) + u^* B_T +  {\g_A-1 \over 2}|u^*|^2T, 
\eeaa
where $u^* := {1+\g_P\over 1+\g_A+\g_P}$ is the corresponding agent's optimal action.

We now consider the  dynamic version of the agent's problem (\ref{PA-VA0}): for  $t\in [0,T]$, 
\bea
\label{PA-VAPt}
 V^A_t(C_T) :=\esssup_{u\in\sU_A} \dbE^{\dbP^u}_t\Big[ U_A\big(C_T  - {1\over 2}\int_t^T |u_s|^2 ds\big)\Big],
\eea
and the principle's problem, given agent's optimal control $u(t, C_T)$:
\bea
\label{PA-VPt}
V^P_t  := \esssup_{C_T\in\sU_P} \dbE^{\dbP^{u(t, C_T)}}_t\Big[U_p(B_T-C_T)\Big],  \q\mbox{subject to}\q V^A_t(C_T) \ge R.
\eea
 Solving the principal's problem (\ref{PA-VPt}) as before
we see that the optimal contract is:
 \beaa
C_T^{t,*} := - {1\over \g_A}\ln (-R) + u^* (B_T-B_t) +  {\g_A-1 \over 2}|u^*|^2(T-t),
\eeaa
where $u^* := {1+\g_P\over 1+\g_A+\g_P}$. Clearly $C_T^{t, *}$ is different from $C^*_T$, thus the 
problem is time-inconsistent.   

Again, the time-inconsistency can be removed if we allow the market value of the agents, the constant $R$, to be time varying
(as it should be!). 
Indeed, if we set
\bea
\label{PA-Rt}
R_t := R\exp\Big(-\g_A[u^*B_t + {\g_A-1\over 2}|u^*|^2 t]\Big),
\eea
and modify the participation constraint of the principal's problem in (\ref{PA-VP0}) to
$V^A_t(C_T) \ge R_t$.
Then the optimal solution to the principle's problem (\ref{PA-VPt}) will become
\beaa
\tilde C_T^{t,*} &=& - {1\over \g_A}\ln (-R_t) + u^* (B_T-B_t) +  {\g_A-1 \over 2}|u^*|^2(T-t)\\
& =& - {1\over \g_A}\ln (-R) + u^* B_T +  {\g_A-1 \over 2}|u^*|^2T = C_T^*.
\eeaa
That is, the problem becomes time-consistent. 

 We note that  the problem \reff{PA-VP0} can also be written as an optimal control problem for a forward-backward SDE.  
To see this, we first note that by some straightforward arguments, one can show that for the optimal contract $C^*_T$, the
identity $V_0(C^*_T) = R$ must hold. 
Therefore we may impose a stronger participation constraint in \reff{PA-VP0}: $V_0(C_T) = R$, and rewrite $Y^A$ as a forward diffusion:
\beaa
Y^A_s = U_A^{-1}(R)   + {\g_A-1\over 2}\int_0^s |Z^A_r|^2  dr + \int_0^s Z^A_r dB_r, \qq s\in[0,T],
\eeaa
which can be thought of as the optimal solution to the agent's problem (\ref{PA-VP0}) with dynamics
\bea
\label{YAu}
Y^{A,u}_s := U_A^{-1}(R)   + {\g_A-1\over 2}\int_0^s |u_r|^2  dr + \int_0^s u_r dB_r, \q s\in[0,T], 
\eea
with the relation $C_T=Y^A_T$. Then, instead of viewing $C_T$ as the principal's control, we may view $u := Z^A$ as the principal's control, and unify the principal-agent problem to the following optimization problem for FBSDEs:
\bea
\label{PA-V0}
V_0 :=\sup_{u\in \sU_A} Y^{\dbP, u}_0,
\eea
where $(Y^{A, u}, Y^{\hP, u})$ is the solution to the (forward) SDE (\ref{YAu}) and the following BSDE
\bea
\label{YPu}
Y^{P,u}_s =  U_P(B_T- Y^{A,u}_T) + \int_s^T u_r Z^{\dbP,u}_r dr -  \int_s^T Z^{\dbP,u}_r dB_r, \qq s\in[0,T],
\eea
respectively.

\begin{rem}
\label{remark2.3}{\rm
The BSDEs appeared in this problems are all one dimensional, thus comparison principle should hold and problem is expected to be time consistent. The time-inconsistency is caused by the fixed constraint
$R=V_0(C_T)$. We removed the time inconsistency by setting $R_t = V^A_t(C^*_T)$ for all $t\in [0, T]$, where $C^*_T= Y^{A.*}_T$ is the optimal contract, which is exactly the  random participation constraint \reff{PA-Rt}.
In more general models, however, the BSDEs could very well be multidimensional, see e.g. \cite{CZ}, and the comparison principle would indeed fail. }
\qed
\end{rem}

\subsection{The probability distortion problem}
\label{sect-distortion}
In this subsection we show that the probability distortion problem considered in \cite{XZ} can also be recast as an optimization problem with controlled BSDEs.  With a slight variation, the problem  in \cite{XZ} can be understood as follows:
\bea
\label{distortion-V0}
V_0 := \sup_{\t} \int_0^\infty  w(\dbP_0(U(B_\t)\ge x)) dx,
\eea
where $\t$ is running over all stopping times, $U \ge 0$ is a utility function, and the probability distortion function $w: [0, 1] \to [0,1]$ is a continuous and strictly increasing function such that $w(0) = 0$ and $w(1)=1$.  If  $w(x) = x$ for all $x\in [0,1]$, then $V_0 = \sup_\t \dbE[U(B_\t)]$, which is a standard optimal stopping problem and is time consistent. However, for general distortion function $w$, the problem is typically time inconsistent as was showed in \cite{XZ}, where the optimal stopping time was constructed by using some quantile  functions and the Skorohod embedding theorem. 

To write 
 \reff{distortion-V0}  in the form of \reff{BSDE-V0}, we let $\t$ be the control and $x\in [0, \infty)$ be the parameter. For each $x$ and $\t$, introduce a BSDE:
\bea
\label{distortion-BSDE}
Y^{x,\t}_t = \1_{\{U(B_\t)\ge x\}} - \int_t^T Z^{x,\t}_sdB_s.
\eea
That is, we view $Y^\t := (Y^{x,\t})_{x\in [0, \infty)}$ as the solution to a (uncountably) infinite dimensional BSDE. Then we have
\bea
\label{distortion-V0phi}
V_0 = \sup_\t \f\big(Y^\t_0\big),\q \mbox{where}\q \f(f)  := \int_0^\infty f(x) dx.
\eea

\subsection{A deterministic example}
\label{sect-deterministic}
It is a common suspicion that the  random uncertainty involved in the underlying problem may play some fundamental role
in  the time inconsistency. To conclude this section we provide a simple deterministic example where the comparison principle fails in order to show that the time inconsistency is more of a structural issue than an information issue.

Let $T> 1$, and $\sU_{[s, t]}$ be the set of deterministic functions $u: [s, t] \to [0,1]$. Consider the deterministic optimization 
problem:
\bea
\label{deterministicVt}
 V_t := \sup_{u\in \sU_{[t, T]}} Y^{1,u}_t,~\mbox{where}~ Y^{1,u}_t := \int_t^T [u_s-Y^{2,u}_s] ds,~ Y^{2,u}_t := \int_t^T u_s ds, ~
t\in [0, T].
\eea
By straightforward calculation, we obtain that
\bea
\label{deterministicY1u}
Y^{1,u}_t = \int_t^T [u_s - \int_s^T u_r dr] ds = \int_t^T [1+t-s] u_s ds,
\eea
and then clearly the optimal control is:
$u^{t, *}_s := \1_{[t, (1+t)\wedge T]}(s)$, $ t\le s\le T$.
In particular, for $0<t<T-1$, we see that 
\bea
\label{deterministicInconsistent}
u^{0,*}_s = 0 \neq 1= u^{t,*}_s,\q  s\in (1, 1+t).
\eea
That is, the problem \reff{deterministicVt} is time inconsistent.

\section{Characterization of Time Consistency in Our Model}
\label{sect-BSDE}
\setcounter{equation}{0}

Having argued in  previous section that many time-inconsistent problems  can be recast as optimization problems  with controlled BSDEs/FBSDEs,
in the rest of the paper we shall focus exclusively on such class of optimization problems and introduce our main schemes. Again, our purpose here is to 
revitalize the ``dynamical programming principle" (DPP) in  a time-inconsistent situation, without assuming the existence of the optimal control. As we pointed out in Introduction, in order to focus more on the main ideas, we shall consider
only the case where the controlled  dynamics are  finite dimensional BSDEs, with the forward component being simply the driving Brownian motion itself.  The extension to controlled forward SDEs requires some heavier notations but no substantial difficulty.
The generalization to infinite dimension is more challenging in general, and we shall leave it to future study.

We begin with a precise description of the framework. Let $U$ be a Polish set, and  $\sU:=\dbL^0_\hF([0,T];U)$.  Consider the following $d'$-dimensional BSDE: 
\bea
\label{Y1}
Y^u_t = \xi + \int_t^T f(s, Y^u_s, Z^u_s, u_s) ds -\int_t^T Z^u_s dB_s, \qq t\in[0,T].
\eea
Now, for a given cost function $\f:\hR^{d'}\to \hR$, we define the following optimization problem:
\bea
\label{V0}
V_0(\xi) := \sup_{u\in \sU} \f(Y^u_0),\q \mbox{for any}~\xi\in \dbL^2(\cF_T;\hR^{d'}).
\eea

Throughout this paper we shall make use of the following {\it Standing Assumptions}:
\begin{assum}
\label{assum-f}
(i) The generator $f: [0, T] \times \O\times \dbR^{d'} \times \dbR^{d'\times d} \times U \to \dbR^{d'}$ is $\dbF$-progressively measurable in all variables, uniformly Lipschitz continuous in $(y,z)$, and 
\beaa
\dbE\Big[\Big(\int_0^T \sup_{u\in U} |f(t,0,0,u)|dt\Big)^2\Big] < \infty.
\eeaa

(ii) The function $\f : \dbR^{d'} \to \dbR$ is continuous.  
\end{assum}

Given $\xi\in \dbL^2(\cF_T;\hR^{d'})$, it is by now well-understood that, under Assumption \ref{assum-f},  BSDE \reff{Y1} is 
well-posed for any   $u\in \sU$,  and $\{Y^u_0, u\in\sU\}$ is a bounded set in $\hR^{d'}$. Thus $V_0(\xi)$ in (\ref{V0}) is well defined. We shall refer to 
problem (\ref{V0}) as the {\it static problem}. 

We now consider the problem (\ref{V0}) in a dynamic setting. For $0\le t\le T$, we define:
\bea
\label{Vt}
V_t(\xi) := \esssup_{u\in \sU} \f(Y^u_t).
\eea
As we observed in the previous section, when $\f$ is non-monotone or when $d'\ge 2$, the problem (\ref{Vt}) is typically 
time inconsistent in the sense that the optimal control of static problem (\ref{V0}) is no longer optimal for the  dynamic problem (\ref{Vt}) over the time duration $[t, T]$. We should note, however, that such a characterization, although self-explanatory and easy to understand, has a
fundamental drawback. That is, it relies on the existence of optimal control, 
which in general is a tall order. In fact, it is by no means clear why problems \reff{V0} and \reff{Vt} will possess any optimal control, 
which in theory would make it impossible to check the time-consistency of the problem.

To get around this deficiency we propose a more generic characterization of time-inconsistency, based on the DPP for the value function. To facilitate our discussion let us introduce another notation. For any $0<t\le T$, $\eta\in \dbL^2(\cF_t)$,  and $u\in \sU$, let $(\sY^u(t, \eta), \sZ^u(t,\eta))$ denote the solution to the following BSDE on $[0, t]$:
\bea
\label{cY}
\sY^u_s = \eta + \int_s^t f(r, \sY^u_r, \sZ^u_r, u_r) dr - \int_s^t \sZ^u_r dB_r,\q 0\le s\le t.
\eea
Clearly, using the notation $\sY^u(\cd, \cd)$ and uniqueness of the solution to BSDE \reff{cY} we can write: $Y^u_s = \sY^u_s(t, Y^u_t)$, $0\le s\le t\le T$; and, in particular, $Y^u_0 = \sY^u_0(t, Y^u_t)$,  $t\in[0,T]$.

\ms

We  illustrate the idea through two examples where  $\f$ is monotone and  the BSDE satisfies the comparison principle.

\begin{eg}
\label{eg-mon1}
{\rm Assume that Assumption \ref{assum-f} is in force, and assume further that $d'=1$ and $\f$ is increasing. 
Then, it is clear that the static problem (\ref{V0}) is equivalent to $V_0(\xi) := \f\Big(\sup_{u\in \sU} Y^u_0\Big)$. 
On the other hand, by the comparison principle of BSDEs and the monotonicity of $\f$, we see immediately that 
the dynamic problem (\ref{Vt}) can also be written as: $V_t(\xi) = \f(\overline Y_t)$, $0\le t\le T$, where 
$\overline f(s,\o, y,z) := \sup_{u\in U} f(s, \o, y,z,u)$, and 
\beaa
 \overline Y_s = \xi + \int_s^T \overline f(s, \overline Y_s, \overline Z_s) ds - \int_s^T \overline Z_sdB_s, \q s\in[0,T].
 \eeaa

We claim that this problem is {\it  time-consistent} in the sense that the following DPP holds: 
\bea
\label{mon1-DPP}
V_{t_1}(\xi) =  \esssup_{u\in \sU} \f( \cY^u_{t_1}(t_2, \overline Y_{t_2})),\q ~ 0\le t_1<t_2\le T.
\eea
Indeed, for simplicity 
we set $t_1:=0$ and $t_2:=t$. For any $u\in \sU$, we write $Y^u_0 = \cY^u_0(t, Y^u_t)$.  By the comparison principle of BSDE and the monotonicity of $\f$, we see that $Y^u_t \le \overline Y_t$ which implies $Y^u_0=\sY^u_0(t, Y^u_t) \le  \sY^u_0(t, \overline Y_t)$ and consequently $\f(Y^u_0) \le \f( \sY^u_0(t, \overline Y_t))$, thanks to the monotonicity of $\f$. 
Since $u$ is arbitrary, we conclude that
\bea
\label{mon1-comparison}
V_0(\xi) \le   \sup_{u\in \sU}\f( \sY^u_0(t, \overline Y_t)).
\eea
To see the opposite inequality of (\ref{mon1-comparison}), for any $\e>0$, we apply the standard measurable selection theorem
to get a measurable function $I_\e: [0, T]\times \O \times \dbR\times \dbR^{1\times d} \to U$ such that
\bea
\label{feqbarf}
f(s, \o, y, z, I_\e(s,\o, y,z)) \ge \overline f(s,\o, y, z) -\e, \q \forall (s,\o, y, z). 
\eea
Set $u^\e_s  := I_\e(s, \overline Y_s, \overline Z_s)$, $t\le s\le T$. By standard BSDE arguments we see that
\bea
\label{1d-Ye}
\ol Y_t \le Y^{u^\e}_t + C\e.
\eea
Now for any $u\in \sU$, by standard BSDE arguments again, it follows from \reff{1d-Ye} that
\beaa
\sY^u_0(t, \overline Y_t) \le \sY^u_0(t,  Y^{u^\e}_t) +C\e= Y^{u\otimes_t u^\e}_0 +C\e\le V_0(\xi) + C\e,
\eeaa
where $u\otimes_t u^\e:= u\1_{[0,t)} + u^\e\1_{[t, T]}$.  By the arbitrariness of $u$ and $\e$, we prove the opposite inequality 
in (\ref{mon1-comparison}), whence the DPP \reff{mon1-DPP}.
\qed}
\end{eg}

We should note that the DPP \reff{mon1-DPP} does not require the existence of optimal control, but it indeed characterizes the time consistency. Moreover, when $U$ is compact  and $f$ is continuous in $u$, there exists a  measurable function 
$I: [0, T]\times \O \times  \dbR\times \dbR^{1\times d} \to U$ such that
\beaa
f(s, \o, y, z, I(s,\o, y,z)) = \overline f(s,\o, y, z), \q \forall (s,\o, y, z). 
\eeaa
In this case, one can easily check that $u^*_s  := I(s, \overline Y_s, \overline Z_s)$ is optimal both for $V_0(\xi)$ and for any $V_t(\xi)$. So the problem is time consistent in terms of optimal control as well.

\begin{rem}
\label{rem-1d-comparison}
{\rm
  As we see in the argument leading to \reff{mon1-comparison}, the DPP \reff{mon1-DPP} clearly relies on the comparison principle of the BSDE and the monotonicity of $\f$. In fact, as we saw in \S\ref{sect-1d}, the comparison alone is not sufficient for the time consistency. 
\qed}
\end{rem}

The next example reinforces the importance of comparison principle for 
time consistency.

\begin{eg}
\label{eg-mon2}
{\rm Let $d' \ge 2$. Consider the following multidimensional BSDE: for $i=1,\cds, d'$,
 \beaa
 \overline Y^i_t = \xi_i + \int_t^T \overline f_i(s, \overline Y_s, \overline Z^i_s) ds -\int_t^T \overline Z^i_s dB_s,
 \eeaa
where $\overline f_i(t,y,z_i) := \sup_{u\in U} f_i(t,y,z_i, u)$. Assume that 

(i) for $i=1,\cds, d'$, $f_i$ does not  depend on $z_j$ and is increasing in $y_j$, for all $j\neq i$; and 

(ii)  $\f$ is increasing in each component. 

\no Then it is well-known that the comparison principle remains true for such BSDEs. Following the similar arguments as in Example \ref{eg-mon1} we can show that $V_t(\xi) = \f(\overline Y_t)$, $0\le t\le T$, and
 \beaa
   V_{t_1}(\xi) =  \esssup_{u\in \sU} \f( \sY^u_{t_1}(t_2, \overline Y_{t_2})),\q \hP\as\neg\neg, \q 0\le t_1<t_2\le T.
  \eeaa
Consequently, the problem is time consistent.
\qed}
\end{eg}

From Example \ref{eg-mon2} we see the essential roles that comparison principle and the monotonicity of some key coefficients
play in the time consistency. 
In general, the comparison principle fails for $d'>2$ except for some special cases. We refer to \cite{HP} for some detailed analysis on this issue. 
We note that the problem will remain time consistent if $f_i$  and $\f$ are monotone on the corresponding variables in a
compatible manner (e.g., $f_i$ is decreasing in $y_j$ and $\f$ is decreasing in all its variables).  The result would be very different if such compatibility is violated. In fact,  as we saw in \S\ref{sect-deterministic}, when $f_i$ is decreasing in $y_j$ but $\f$ is increasing, the problem becomes time inconsistent. 
%

To study the general time-inconsistent problem we
propose the following definition. 
\begin{defn}
\label{defn-consistent}
An $\dbF$-progressively measurable function  $\Phi: [0, T]\times \O\times \dbR^{d'}\to \dbR$ is called a ``time consistent dynamic utility function" for problem \reff{Y1}-\reff{V0} if 

(i)  $\Phi(0, y) = \f(y)$, 

(ii) there exists  a mapping $\ol Y:[0,T]\times \O\mapsto \hR^{d'}$ satisfying $\ol Y_t\in \hL^2(\cF_t; \dbR^{d'})$, for $t\in[0,T]$ and $\ol Y_T = \xi$, 
$\hP$-a.s., such that the following DPP holds:
 \bea
 \label{olYDPP}
\Phi(t_1, \ol Y_{t_1})  =  \esssup_{u\in \sU} \Phi(t_1,  \sY^u_{t_1}(t_2, \overline Y_{t_2})),  \qq 0\le t_1<t_2\le T.
  \eea
  
  \no In particular, in this case we  say that the following dynamic processes is time consistent:
  \bea
\label{tildeVt}
\tilde V_t(\xi) := \Phi(t, \ol Y_t) =  \esssup_{u\in \sU} \Phi(t,Y^u_t).
\eea
  \end{defn}

 \begin{rem}
 \label{rem-forwardutiltiy}
 {\rm 
The time consistent dynamic utility function $\Phi$ is motivated in part by the notion of the forward utility proposed in \cite{MZ1, MZ2, EM}, because both evolve forwardly in time. 
It should be noted, however, that there is a fundamental difference here: for each $t\in [0, T]$, the forward utility $U(t,\cd)$  in \cite{MZ1, MZ2, EM} acts on $t$ and optimizes
over the time duration $[0,t]$, whereas our dynamic utility $\Phi(t,\cd)$ acts on terminal time $T$ and optimizes over the
time duration $[t, T]$. 
\qed}
\end{rem}

We would like to emphasize the following three  three main features of Definition \ref{defn-consistent}:

1) $\tilde V_0(\xi) = V_0(\xi)$, thanks to condition (i). This means the dynamic problem is consistent with the static problem.

2)  The function $\Phi$ is defined ``forwardly", with an initial value, and the mapping $\ol Y$ is defined backwardly, with a terminal value. We should particularly note that at this point we \underline{\it do not} require the $t$-measurability of the mapping $\ol Y$; and

3) The time consistency is characterized by the DPP, which does not require the existence of optimal control.   

It is easy to see that the function $\Phi(t,\cd) \equiv\f$ in  Examples \ref{eg-mon1} and \ref{eg-mon2}  is a time consistent dynamic utility. Furthermore, if the optimal control $u^*$ exists,  we may simply set $\ol Y := Y^{u^*}$, and in this case one can easily find a desired $\Phi$, as we see in the examples in previous section. However, in general, we need to find the $\ol Y$ whose dynamics (if it exists)  may help us to either determine the optimal control $u^*$, if any, or find conditions for the existence of optimal control.
We should also note that
 the dynamic utility function $\Phi$ is not unique. In fact, if $\Phi$ is a time consistent dynamic utility, then for any process $\th$ with $\th_0=0$, $\tilde \Phi(t, y) := \Phi(t, y) + \th_t$ is also a time consistent dynamic utility. Since our main difficulty is the existence of such $\Phi$, in this paper we impose minimum requirements on $\Phi$. 

%

In the rest of this paper, we shall propose three possible approaches to attack the general time inconsistent optimization problems (in the sense that $\Phi(t,\cd) \equiv\f$ is not a time consistent dynamic utility function). Each approach has its pros and cons. We  note that in this
paper we focus mainly on the ideas, rather than the actual solvability of the resulting problems, which could be highly technical, and
 may call for some new developments in the respective areas.

\section{The Duality  Approach}
\label{sect-duality}
\setcounter{equation}{0}

\subsection{Heuristic analysis in Markovian case}
In this section we present a duality approach that is simple but quite effective if one focuses only on finding the value of the static problem
(\ref{V0}). To illustrate the idea better we begin by  considering the Markovian case, that is, we assume that in BSDE (\ref{Y1}) $\xi = g(B_T)$ and $f = f(t, B_t, y,z,u)$. We shall start with a heuristic arguments, and give the proof for the general non-Markovian (or say path-dependent) case.

To begin with, for each $(t,x)\in[0,T]\times\hR^d$, 
consider the set
\bea
\label{cDtx}
 \sD(t,x) := \Big\{y\in \dbR^{d'}: \exists Z \in \dbL^2_\hF([0,T]), u\in\sU_{[t, T]}, ~\mbox{s.t.~} X^{t,x,y,Z,u}_T = g(B^{t,x}_T),~\dbP_0\mbox{-a.s.}\Big\},
\eea
where $B^{t,x}_s:= x+B_s-B_t$, $s\ge t$, and $X^{t,x,y,Z,u}$ is the solution to the forward SDE:
\bea
\label{X}
X_s 
&=& y - \int_t^s f(r, B^{t,x}_r,  X_r, Z_r, u_r) dr + \int_t^s Z_r dB_r, \qq t\le s\le T.
\eea
Clearly, 
$X$ can be thought of as a  forward version of the solution to the BSDE (\ref{Y1}) on $[t, T]$, and the set $\sD(t,x)$ is simply 
the {\it reacheable set} $\{Y^u_t, u\in \sU\}$ given $B_t=x$.  In particular, $\sD(0,0) = \{Y^u_0: u\in \sU\}$, 
and  our original optimization can be written as 
\bea
\label{cDtxV0}
 V_0(\xi) = \sup_{y\in \sD(0,0)} \f(y).
\eea
It is worth noting that $ \sup_{y\in \sD(0,0)} \f(y)$ in \reff{cDtxV0} is a finite dimensional optimization problem. So the 
value $V_0(\xi)$ could be determined rather easily, provided one can characterize the set $\sD(0,0)$, which we now describe. 

To this end, we borrow the idea of the {\it method of optimal control} for solving a forward-backward SDE
(cf. \cite{MY}). Consider the following dual control problem:
\bea
\label{Wtxy}
W(t,x,y) := \inf_{Z, u} \dbE\Big\{\big|X^{t,x,y,Z,u}_T - g(B^{t,x}_T)\big|^2\Big\}.
\eea
 Clearly, (\ref{Wtxy}) is a standard stochastic control problem, and it is well-known that $W$ should be the (unique) viscosity solution to the following (degenerate) HJB equation:
\bea
\label{WPDE}
\left\{\ba{lll}
\dis \pa_t W + {1\over 2} \pa^2_{xx} W + \inf_{z, u}\Big\{{1\over 2} \pa^2_{yy}W : (z z^\top) + \pa^2_{xy}W : z - \pa_y W \cd f(t,x,y,z,u)\Big\}=0;\ms\\
\dis W(T,x,y) = |y-g(x)|^2.
\ea\right.
\eea
By definition (\ref{cDtx}) it is clear that $W(t,x,y)=0$ whenever $y\in \sD(t,x)$. More generally, we expect and will show that, for any $(t,x)$,
the following duality relationship between  the set $\sD(t,x)$ and the ``nodal set" of the function $W$ holds:
\bea
\label{WDtx}
\sN(t,x) := \Big\{y\in \dbR^{d'}: W(t,x,y) = 0\Big\} ~=~ \overline{\sD(t,x)}.
\eea
where $\overline{\sD(t,x)}$ denotes the closure of $\sD(t,x)$.  Then \reff{cDtxV0} amounts to saying that
\bea
\label{V0Ntx}
V_0(\xi) = \sup_{y\in \sD(0,0)} \f(y)=\sup_{y\in \sN(0,0)} \f(y).
\eea
In other words, we have characterized the set $\sD(0,0)$ in terms of $\sN(0,0)$, the nodal set of $W$, which is 
a much benign task to deal with (for example, numerically).  Moreover, note that the nodal set $\sN(0,0) \subset \dbR^{d'}$ is closed, then the above optimization problem  has a maximum argument $y^*\in \sN(0,0)$. Consequently, the static optimization problem \reff{V0} has an optimal control if and only if there exists $y^*\in \sD(0,0)$.

\begin{rem}
\label{remark4.1}{\rm
(i) An important ingredient in the duality approach is 
the ``reachable set" $\sD(\cd, \cd)$. Unlike the standard optimal control literature where 
reachable sets are temporally forward, it is easy to see from (\ref{cDtx}) that the family  $\{\sD(t,\cd)\}_{0\le t\le T}$ is a backward, set-valued dynamic system with terminal condition $\sD(T,x) = \{g(x)\}$, and  
as we shall see later in this section, it satisfies  a  {\it geometric DPP} in  the spirit of \cite{ST}.

(ii) The duality approach could be combined with the time consistency 
in the sense of Definition \ref{defn-consistent} 
as follows.  Assuming we could find a desired time consistent dynamic utility $\Phi$, which we hope will take the form $\Phi(t, B_t, y)$ in the Markovian case, then by the duality \reff{WDtx} we have the following time consistent value function:
\bea
\label{dual-tildeV}
\tilde V_t(\xi) = \esssup_{y\in \sN(t, B_t)}  \Phi(t, B_t, y).
\eea
Moreover, since the nodal set $\sN(t, B_t) \subset \dbR^{d'}$ is closed,  the above optimization problem  has maximum argument $\ol Y_t$, which would serve for the purpose of Definition \ref{defn-consistent}.  The task of finding a desired $\Phi$, however, is challenging. We shall investigate it in next section.
\qed}
\end{rem}


\subsection{The duality approach for the general path dependent case}
We now carry out the  duality approach rigorously in the general  path dependent (or  non-Markovian) case. To begin with, we recall the canonical set-up  introduced in the beginning of \S\ref{sect-inconsistent}. Moreover, for any $t\in [0, T]$, denote by $\O^t:= \{\o\in C([t, T], \dbR^d): \o_t=0\}$ the shifted canonical space on $[t, T]$, and define $B^t, \dbF^t, \dbP_0^t, \L^t$,  $\sU^t$ etc on $\O^t$ in obvious sense. 
Furthermore, for any $\o\in \O$ and $\tilde \o\in \O^t$, we introduce the concatenation: $\hat\o:=\o\otimes_t\tilde\o=\o{\bf 1}_{[0,t]}+(\o_t + \tilde \o){\bf 1}_{[t, T]}$. 
Moreover, for $\xi\in \hL^0(\O)$ and $(t,\o) \in \L$, denote $\xi^{t, \o}(\tilde\o):=\xi(\o\otimes_t \tilde\o)$, for all $\tilde\o\in\O^t$.

Similar to \reff{cDtx},  for any $(t,\o) \in \L$ we define, 
\bea
\label{cDto}
\sD(t,\o) := \Big\{y\in \dbR^{d'}: \exists (Z, u) \in \dbL^2(\dbF^t, \dbR^{d'\times d}) \times \sU^t,~\mbox{s.t.~} X^{t,\o,y,Z,u}_T = \xi^{t,\o},~\dbP_0^{t}\mbox{-a.s.}\Big\},
\eea
where $X^{t,\o,y,Z,u}$ is the solution to the following (forward) SDE:
\bea
\label{Xo}
X_s = y - \int_t^s f^{t,\o}(r,  B^t_\cd,  X_r, Z_r, u_r) dr + \int_t^s Z_r dB^t_r, \q t\le s\le T,~\dbP^{t}_0\mbox{-a.s.}
\eea
Here the function $f^{t, \o}(r, \tilde\o, y, z, u)$, $(r, \tilde\o)\in\L^t$ is defined the same as $\xi^{t,\o}$ 
explained before. 
Again, it is easy to see that $\sD(0,0) = \{Y^u_0: u\in \sU\}$ remains true. Thus we still have
\bea
\label{cDtoV0}
V_0(\xi) = \sup_{y\in \sD(0,0)} \f(y).
\eea
We now introduce a dual control problem in the path-dependent setting:
\bea
\label{Wtoy}
W(t,\o,y) := \inf_{(Z, u)\in  \dbL^2(\dbF^t, \dbR^{d'\times d}) \times \sU^t} \hE^{\hP_0^{t}}\Big[\big|X^{t,\o,y,Z,u}_T - \xi^{t,\o}\big|^2\Big].
\eea
 Our main duality result is as follows. 
\begin{thm}
\label{thm-duality}
Let Assumption \ref{assum-f} hold, and assume further that,  for any $(t,\o)\in \L$,
\bea
\label{duality-assum}
\dbE^{\dbP_0^t}\Big[\Big(\int_t^T \sup_{u\in U} |f^{t,\o}(s, B^t_\cd, 0,0, u)|ds\Big)^2 + |\xi^{t,\o}|^2\Big] < \infty.
\eea
 Then,  for any $(t,\o)\in \L$, we have
\bea
\label{WDto}
\sN(t,\o) := \Big\{y\in \dbR^{d'}: W(t,\o,y) = 0\Big\} ~=~ \overline{\sD(t,\o)}.
\eea
Consequently,
$V_0(\xi) = \sup_{y\in \sN(0,0)} \f(y)$.
\end{thm}

\proof Noting  \reff{cDtoV0}
and the continuity of $\f$, we shall prove only \reff{WDto}. 

We first prove the regularity of $W$ in $y$: for any $(t,\o)\in\L$, and $y_1, y_2\in\hR$,
\bea
\label{Wyreg}
 |W(t,\o, y_1)-W(t,\o,y_2)|\le C(t,\o) [1+|y_1|+|y_2|] |y_1-y_2|,
\eea
where $C(t,\o)>0$ is independent of $y$.
Indeed, by (\ref{Wtoy}) and (\ref{duality-assum}), it is readily seen that
\beaa
W(t,\o, y) \le C(t,\o) [1+|y|^2].
\eeaa
Now for any $0<\e<1$, we choose $(Z^\e, u^\e)\in  \dbL^2(\dbF^t, \dbR^{d'\times d}) \times \sU^t$ such that 
\beaa
\hE^{\hP_0^{t}}\Big[\big|X^{t,\o,y_2,Z^\e,u^\e}_T - \xi^{t,\o}\big|^2\Big] \le W(t,\o, y_2) + \e \le C(t,\o)[1+|y_2|^2].
\eeaa
 By the standard BSDE arguments, it is then clear that, under Assumptions \ref{assum-f}, we have
\beaa
\hE^{\hP_0^{t}}\Big[\big|X^{t,\o,y_1,Z^\e,u^\e}_T -X^{t,\o,y_2,Z^\e,u^\e}_T \big|^2\Big] \le C|y_1-y_2|^2. 
\eeaa
Then, denoting $X^i:= X^{t,\o,y_i,Z^\e,u^\e}$, $i=1,2$, we have
\beaa
W(t,\o, y_1)-W(t,\o,y_2) &\le&  \hE^{\hP_0^{t}}\Big[\big|X^1_T - \xi^{t,\o}\big|^2\Big] - \hE^{\hP_0^{t}}\Big[\big|X^2_T - \xi^{t,\o}\big|^2\Big] + \e\\
&\le&\hE^{\hP_0^{t}}\Big[\big|X^1_T -X^2_T \big|^2 +2|X^1_T-X^2_T| |X^2_T-\xi^{t,\o}|\Big] + \e\\
&\le& C|y_1-y_2|^2 + C(t,\o) [1+|y_2|]|y_1-y_2|+\e \\
&\le& C(t,\o)[1+|y_1|+|y_2|]|y_1-y_2| +\e.
\eeaa
Since $\e$ is arbitrary, we obtain the desired estimate (\ref{Wyreg}) for $W(t,\o, y_1)-W(t,\o,y_2)$. Switching the roles of 
$y_1$ and $y_2$ we can also obtain the estimate for $W(t,\o, y_2)-W(t,\o,y_1)$, whence \reff{Wyreg}.

Next, we  fix $(t, \o)\in\L$ and let $y\in \sD(t,\o)$. By definition there exists $(Z, u)\in \dbL^2(\dbF^t, \dbR^{d'\times d}) \times \sU^t$ such that $X^{t,\o, y, Z, u}_T = \xi^{t,\o}$, $\dbP^{t}_0$-a.s. Then we must have
\beaa
W(t,\o,y) \le \hE^{\dbP^{t}_0}\Big[\big|X^{t,\o, y, Z, u}_T - \xi^{t,\o}\big|^2\Big] = 0.
\eeaa
That is, $y\in \sN(t,\o)$ and consequently $\sD(t,\o) \subset \sN(t,\o)$. Moreover, the $y$-continuity of $W$ in  \reff{Wyreg}  then implies that
$\sN(t,\o)$ is a closed set, which leads to that $\overline{\sD(t,\o)} \subset  \sN(t,\o)$.

Conversely, if $y\in \sN(t,\o)$, then by definition for any $\e>0$, there exists $(Z^\e, u^\e)\in \dbL^2(\dbF^t, \dbR^{d'\times d}) \times \sU^t$ such that 
\bea
\label{Zue}
\hE^{\dbP^{t}_0}\Big[\big|\xi^{t,\o}_\e - \xi^{t, \o}\big|^2\Big] \le \e,\q \mbox{where}\q \xi^{t,\o}_\e := X^{t,\o,y,Z^\e,u^\e}_T.
\eea
Now by the standard BSDE estimates we have, for the given $(t,\o)\in\L$,
\beaa
|Y^{u^\e}_t(\o) - y|^2 = \Big|\sY^{u^\e}_t(T, \xi^{t,\o}) - \sY^{u^\e}_t(T, \xi^{t,\o}_\e)\Big|^2 \le C \hE^{\dbP^{t}_0}\Big[\big|\xi^{t,\o}_\e - \xi^{t,\o}\big|^2\Big] \le C\e.  
\eeaa
Since $Y^{u^\e}_t(\o)\in \sD(t,\o)$ and $\e$ is arbitrary, we see that $y\in \overline{\sD(t,\o)}$.
\qed


\subsection{Characterization of $W$ by PPDEs}
\label{sect-path}
It is well understood that, in Markovian case, the dual value function $W$ is the viscosity solution to HJB equation \reff{WPDE}. In this subsection we extend this characterization of $W$ to path dependent case via the newly established viscosity theory developed in  
\cite{EKTZ, ETZ1,ETZ2}.
The path derivatives introduced here will also be important in \S\ref{sect-master}. Since the results here are irrelevant to the rest of the paper, 
we shall focus only on the main ideas without getting into all the technical details. The interested reader is referred to \cite{ETZ1, ETZ2} for more on pathwise analysis involved in the arguments.

We first consider the following pseudo-metric on $\O$ and $\L$ introduced in  
\cite{Dupire} and \cite{CF}:
 \bea
 \label{dinfty}
\|\o\|_t := \sup_{0\le s\le t}|\o_s|,\q  d_\infty((t,\o), (t', \o')) := |t-t'|^{1\over 2} + \|\o_{t\wedge \cd} - \o'_{t'\wedge \cd}\|_T.
\eea
Let $C^0(\L)$ be the set of processes $v: \L\to \dbR$ that are continuous under $d_\infty$.
We note that any $v\in C^0(\L)$ is $\dbF$-progressively measurable. When $v$ is taking values in, say, $\hR^k$, we denote it by $C^0(\L;\hR^k)$. 
Let $\dbS^d$ denote the set of $d\times d$-symmetric matrices. We say a probability measure $\dbP$ on $\O$ is a semi-martingale measure if $B$ is a semimartingale under $\dbP$.  We now introduce the path derivatives for processes, which is due to \cite{ETZ1} and inspired by \cite{Dupire}.
 \begin{defn}
\label{defn-C12}
Let $v\in C^0(\L)$. We say $v\in C^{1,2}(\L)$ if there exist $\pa_t v\in C^0(\L; \hR), \pa_\o v\in C^0(\L; \hR^d), \pa^2_{\o\o} v \in C^0(\L; \hS^d)$ such that the following functional Ito formula holds: for any semimartingale measure $\hP$,
\bea
\label{FIto}
d v(t,\o) = \pa_t v dt + \pa_\o v \cd dB_t + {1\over 2} \pa^2_{\o\o}v : d\la B\ra_t,\q\hP\mbox{-a.s.}
\eea
 \end{defn}
\no We remark that the path derivatives $\pa_t v, \pa_\o v, \pa^2_{\o\o} v$, if they exist, are unique.

Notice that the function $W$ in \reff{Wtoy} is defined on $\L \times \hR^{d'}$. By increasing the space dimension and viewing $y$ as the current value of the additional paths, one may easily extend all the above notions for functions on $\L\times \hR^{d'}$
(see \cite{ETZ1} for details).

We shall make use of the following extra assumption:
\begin{assum}
\label{assum-ocont}
(i) The mapping $(t, \o)\mapsto f(t, \o,  y,z, u)$  is uniformly continuous under $d_\infty$, uniformly in $(y,z,u)$, and $f(t,\o, 0,0, u)$ is bounded;

(ii) The mapping  $\o\mapsto \xi(\o)$  is uniformly continuous under $\|\cd\|_T$ and is bounded.
\end{assum}

Under Assumption \ref{assum-ocont}, by standard BSDE arguments one can easily show that the function $W$ defined by \reff{Wtoy} is uniformly continuous and bounded. 
 It then follows from \cite{ETZ1} that $W$ is a viscosity solution of the following path dependent HJB equation: 
\bea
\label{WPPDE}
\left\{\ba{lll}
\dis 0= \pa_t W + {1\over 2}\tr( \pa^2_{\o\o} W)
+ \inf_{(z, u)}\Big[{1\over 2} \pa^2_{yy}W\neg :\neg (z z^\top) + \pa^2_{\o y}W\neg \cd\neg z - \pa_y W \neg\cd\neg f(t,\o, y,z,u)\Big];\ms\\
W(T,\o,y) = |y-\xi(\o)|^2.
\ea\right.
\eea
In particular, if $W\in C^{1,2}(\L\times \dbR^{d'})$, then $W$ is a classical solution to the above PPDE.

We shall remark though, the above PPDE is degenerate, and thus the uniqueness result of \cite{ETZ2} does not apply here. We refer to the more recent works \cite{RTZ,  EZ}, in which it was shown that $W$ is indeed the unique viscosity solution. We also refer to \cite{RT, ZZ}
for numerical methods for PPDEs.

We conclude this section by providing a rigorous form of the ``geometric DPP" for the set valued process $\sD(t,\o)$ defined by \reff{cDto}, that has been instrumental in the discussions of this section.  Intuitively, in light of \cite{ST}, we expect the following identity: 
\bea
\label{cDtoDPP}
\sD(t_1,\o) &=& \Big\{y\in \hR^{d'}: \exists (Z, u)\in \dbL^2(\hF^{t_1}, \hR^{d'\times d}) \times \sU^{t_1}~\mbox{such that} \\
 &&\qq\qq X^{t_1,\o,y,Z,u}_{t_2}\in \sD(t_2, \o\otimes_t  B^{t_1}),~\dbP^{t_1}_0\mbox{-a.s.}\Big\},\q 0\le t_1 < t_2 \le T.\nonumber
\eea
Denoting   the right side of \reff{cDtoDPP} by $\sD'(t_1,\o)$, one can easily prove that $\sD(t_1, \o) \subset  \sD'(t_1,\o)$. 
However,  the opposite inclusion is far from obvious. In what follows we prove a weaker version of geometric DPP.
We first recall \reff{WDto} and define, for any $\e>0$,
\bea
\label{cNe}
\sN_\e(t,\o) := \{y\in \dbR^{d'}: W(t,\o, y) \le \e\}.
\eea
It is clear that 
$\sN(t,\o) = \cap_{\e>0} \sN_\e(t,\o)$.

\begin{thm}
\label{thm-GDPP} 
Under Assumptions \ref{assum-f} and \ref{assum-ocont}, the following geometric DPP holds true:  
\bea
\label{cNtoDPP}
\sN(t_1,\o) &=& \bigcap_{\e>0}\Big\{y\in \hR^{d'}: \exists (Z^\e, u^\e)\in \hL^2(\dbF^{t_1}, \dbR^{d'\times d}) \times \sU^{t_1}~\mbox{such that} \\
 &&\q X^{t_1,\o,y,Z^\e,u^\e}_{t_2}(\tilde \o)\in \sN_\e(t_2, \o\otimes_t  \tilde \o),~\hP^{t_1}_0\mbox{-a.e.}~\tilde \o\in \O^{t_1}\Big\},\q 0\le t_1 < t_2 \le T.\nonumber
\eea
\end{thm}
\proof For simplicity, we assume $t_1=0$ and $t_2 = t$, and let $\sN'(t_1,\o)$ denote the right side of \reff{cNtoDPP}. Noting that $\o_0=0$, we shall prove that 
\bea
\label{cN00}
\sN(0,0) = \sN'(0,0) &:=& \bigcap_{\e>0}\Big\{y\in \hR^{d'}: \exists (Z^\e, u^\e)\in \dbL^2(\hF, \hR^{d'\times d}) \times \sU~\mbox{such that} \nonumber\\
&&\qq X^{0,0,y,Z^\e,u^\e}_t(\o) \in \sN_\e(t, \o),~\hP_0\mbox{-a.e.}~ \o\in \O\Big\}. 
\eea
Following the arguments in \cite{ETZ1}, one shows that  $W$ is uniformly continuous in $(t,\o, y)$ with modulus of continuity function $\rho_W(\cd)$, and satisfies the following DPP:
\bea
\label{WDPP}
W(0,0,y) = \inf_{(Z, u)\in \hL^2(\hF, \hR^{d'\times d}) \times \sU} \dbE^{\dbP_0}\Big[W(t, B_\cd, X^{0,0,y, Z, u}_t)\Big].
\eea

Now let $y\in \sN'(0,0)$. For any $\e>0$, let $(Z^\e, u^\e)$ be as in the right side of \reff{cN00}.  Then $W(t, B, X^{0,0,y,Z^\e,u^\e}_t) \le \e$ $\dbP_0$-a.s. and thus $\hE^{\dbP_0}\big[W(t, B, X^{0,0,y,Z^\e,u^\e}_t)\big] \le \e$. This, together with \reff{WDPP}, implies that $W(0,0,y) =0$. Then $y\in \cN(0,0)$ and hence $\sN'(0,0) \subset \sN(0,0)$.

 To see the opposite inclusion, let $y\in \sN(0,0)$, and for any $\e>0$, choose $y_\e\in \sD(0,0)$, such that $|y_\e - y|\le \e$.  By \reff{cDto}, let $(Z^\e, u^\e) \in  \dbL^2(\hF, \hR^{d'\times d}) \times \sU$ be such that $X^{0,\e}_T:=X^{0,0,y_\e,Z^\e,u^\e}_T = \xi$, $\dbP_0$-a.s. It is straightforward to see that, for $\dbP_0$-a.e. $\o\in \O$ and $t\in[0,T]$, 
$(Z^{\e, t,\o},u^{\e, t,\o}) \in \dbL^2(\dbF^{t}, \hR^{d'\times d})\times \sU^t$, and $ (X^{0,\e}_s)^{t,\o} = X^{t, \o , X^{0,\e}_t, Z^{\e,t,\o},u^{\e, t,\o}}_s$, $t\le s\le T$, $\dbP_0^{t}$-a.s.
Consequently, we have $X^{t, \o , X^{0,\e}_t, Z^{\e,t,\o},u^{\e, t,\o}}_T = (X^{0,\e}_T)^{t,\o}   = \xi^{t,\o}$ and thus $X^{0,\e}_t \in \sD(t,\o)$. Now denote $X^{0,y,\e}:=X^{0,0,y,X^\e, u^\e}$, and let $\D X:= X^{0,\e} - X^{0,y,\e}$. Then
\beaa
\D X_s = y_\e- y + \int_0^t \a_r \D X_r dr, \q 0\le s\le t,
\eeaa
where $\a$ is a bounded $\hF$-adapted process, thanks to the Lipschitz continuity of $f$ in $y$. Then clearly $|\D X_t|\le C|y_\e-y|\le C\e$, and thus 
\beaa
|W(t, \o, X^{0,y,\e}_t(\o))| = |W(t, \o, X^{0,y,\e}_t(\o)) - W(t,\o, X^{0,\e}_t(\o))| \le \rho_W(|\D X_t(\o)|) \le \rho_W(C\e).
\eeaa
This implies that $X^{0,0,y,Z^\e,u^\e}_t(\o)\in \sN_{\rho(C\e)}(t,\o)$. Since $\e>0$ is arbitrary, we obtain $y \in \sN'(0,0)$, and thus $\sN(0,0) \subset \sN'(0,0)$.
\qed

\section{The Dynamic Utility  Approach}
\label{sect-forward}
\setcounter{equation}{0}

As we have pointed out in the Introduction, as well as in Definition \ref{defn-consistent}, one of the essential points
in our scheme is to determine the ``time consistent dynamic utility" $\Phi$. We devote this section to the discussion of
its existence.
%
\subsection{The deterministic case}
\label{sect-deterministic2}
We begin with the case where both $f$ and $\xi$ are deterministic, and the admissible controls are also deterministic measurable functions $u\in \hL^0([0, T];U)$.  We shall still assume  Assumption \ref{assum-f} holds, and try
to construct  $\Phi$ explicitly.

Since $\xi$ is deterministic, for $u\in \hL^0([0,T];U)$, the solution to the BSDE (\ref{Y1}), $(Y^u, Z^u)$,
must satisfy $Z^u\equiv0$. Further, if we consider the (deterministic) optimization problem: 
\bea
\label{deterministic-Phi}
\Phi(t, y) := \sup_u \f(Y^{t,y,u}_0), \q\mbox{where \, $Y^{t,y,u}_s = y + \int_s^t f(r, Y^{t,y,u}_r, 0, u_r)dr,\qq 0\le s\le t$}
\eea
then $\Phi$ will be time consistent in the sense that it satisfies the DPP: 
\bea
\label{deterministic-DPP}
\Phi(t_2, y) := \sup_u \Phi(t_1, Y^{t_2,y,u}_{t_1}), \qq \mbox{for $0\le t_1<t_2\le T$}.
\eea
We shall argue  that $\Phi$ is  a time consistent dynamic utility in the sense of Definition \ref{defn-consistent}, 
by identifying the required mapping $\ol Y$. Indeed, note that $\Phi(T,\xi) = V_0(\xi) = \sup_u \f(Y^u_0)$, there exists $u^\e$ such that $\lim_{\e\to 0} \f(Y^{u^\e}_0) = \Phi(T,\xi)$. Denote $\ol f_t := \sup_{u\in U} |f(t,0,0,u)|$. By Assumption \ref{assum-f} we see that $\int_0^T \ol f_t dt <\infty$. One may easily check that
\beaa
\sup_\e \sup_{0\le t\le T}|Y^{u^\e}_t| \le C, \q \sup_\e |Y^{u^\e}_t - Y^{u^\e}_s| \le C\int_s^t [\ol f_r + 1] dr, ~0\le s<t\le T.
\eeaa
Now, applying the Arzela-Ascoli theorem we have, possibly along a subsequence (still denoted by $u^\e$), $\lim_{\e\to 0} \sup_{0\le t\le T} |Y^{u^\e}_t-\ol Y_t|=0$, and $\ol Y$ is an absolutely continuous function.
 
It is clear that $\Phi(0, y) = \f(y)$ and $\ol Y_T = \xi$. Further, for any two functions $u^1, u^2$, denote $u^1\otimes_t u^2 := u^1\1_{[0, t)} + u^2\1_{[t, T]}$. By stability of ODEs, one can easily check that 
\beaa
\Phi(t, \ol Y_t) = \sup_u \f(Y^{t, \ol Y_t, u}_0) = \lim_{\e\to 0} \sup_u \f(Y^{u\otimes_t u^\e}_0).
\eeaa
Now on one hand, we have $ \f(Y^{u\otimes_t u^\e}_0) \le V_0(\xi)$ for any $u$ and $\e$. But on the other hand,  
\beaa
\lim_{\e\to 0} \sup_u \f(Y^{u\otimes_t u^\e}_0)\ge  \lim_{\e\to 0} \f(Y^{u^\e\otimes_t u^\e}_0) =  \lim_{\e\to 0} \f(Y^{u^\e}_0) = V_0(\xi).
\eeaa
Namely, $\Phi(t, \ol Y_t) = V_0(\xi)$. For $0\le t_1 < t_2 \le T$, we can follow the similar arguments to get
\beaa
\sup_u \Phi(t_1, \sY^u_{t_1}(t_2, \ol Y_{t_2})) &=& \lim_{\e\to 0} \sup_u\Phi(t_1, \sY^u_{t_1}(t_2, Y^{u^\e}_{t_2})) =  \lim_{\e\to 0} \sup_u\Phi(t_1, Y^{u\otimes_{t_2} u^\e}_{t_1}) \\
&=&   \lim_{\e\to 0} \sup_u \sup_{u'}\f(Y^{u'\otimes_{t_1} u\otimes_{t_2} u^\e}_0) = V_0(\xi) = \Phi(t_1, \ol Y_{t_1}).
\eeaa
This verifies \reff{olYDPP}. To wit,  $\Phi$ is indeed a time consistent dynamic utility.
\qed

 \begin{rem}
 \label{rem-deterministic}
 {\rm If we denote $\tilde \Phi(t, y) := \Phi(T-t, y)$ and $\tilde f(t,y,z,u) := f(T-t, y,z,u)$, then
 \beaa
 \tilde \Phi(t, y) = \sup_u \f(X^{t, y, u}_T),\q\mbox{where}\q X^{t,y,u}_s = y - \int_t^s \tilde f(r, X^{t,y,u}_r, 0, u_r)dr,\q t\le s\le T.
\eeaa
This is a very standard (deterministic) control problem on $[0, T]$ with utility function $\f$.
However, such a ``time change"  technique would fail in the stochastic case (e.g., when $\xi$ is random), due to the adaptedness requirement.  The master equation approach in \S\ref{sect-master} will address this issue.
\qed}
\end{rem}

\subsection{Dynamic utility via comparison principle}

As we saw in \S\ref{sect-BSDE}, especially  Examples \ref{eg-mon1} and \ref{eg-mon2},  the comparison principle plays a crucial role for time consistency. In this subsection we explore the impact of the comparison principle to the existence of the 
time consistent dynamic utility $\Phi$. To this end, we propose the following slightly stronger form of 
comparison principle:
\begin{defn}
\label{defn-comparison}
We say a mapping $\Phi: \L \times \dbR^{d'}\to \dbR$ satisfies the comparison principle if for any $t_1 < t_2$ and any $\eta, \tilde \eta\in \dbL^2(\cF_{t_2})$, $\Phi(t_2, \eta) \le \Phi(t_2,  \tilde \eta)$, 
$\dbP_0$-a.s. implies that 
\bea
\label{Phicomparison}
 \esssup_{u\in\sU} \Phi(t_1, \sY^u_{t_1}(t_2, \eta)) \le \esssup_{u\in\sU} \Phi(t_1, \sY^u_{t_1}(t_2, \tilde\eta)), \q \dbP_0\mbox{-a.s.}
\eea
\end{defn}
The main result of this subsection is the following theorem.
\begin{thm}
\label{thm-comparison}
Let Assumptions \ref{assum-f} and \ref{assum-ocont} hold and assume there exists a random field $\Phi$ satisfying the following properties:

(i) the mapping $y\mapsto \Phi(t,\o,y)$ is continuous, for fixed $(t, \o)\in[0,T]\times \O$;

(ii) $\Phi(0, \cd, y) = \f(y)$, $\hP_0$-a.s.; and

(iii)  $\Phi$ satisfies the comparison principle \reff{Phicomparison}. 

\no Then $\Phi$ is a time consistent dynamic utility in the sense of Definition \ref{defn-consistent}.
\end{thm}
\proof We shall follow the similar ideas used for the duality approach in previous section, but here we will focus more on the
measurability issue. To this end we 
adjust the notations slightly. For any $t\in [0, T]$, $\eta \in \dbL^2(\cF_t, \dbR^{d'})$, $Z\in \dbL^2(\dbF, \dbR^{d'\times d})$, and 
$u\in \sU$, we denote $X^{t,\eta, Z, u}$ to be the solution to the following random differential equation: 
\bea
\label{XODE}
X^{t,\eta, Z, u}_s = \eta - \int_t^s f(r, X^{t,\eta, Z, u}_r, Z_r, u_r) dr + \int_t^s Z_r dB_r, \q t\le s\le T,~\dbP_0\mbox{-a.s.}
\eea
Clearly, (\ref{XODE}) is essentially an ODE, which can be solved $\o$-wisely. Now define
\beaa
\tilde W(t, y) := \essinf_{(Z,u)\in  \hL^2(\hF, \dbR^{d'\times d}) \times \sU} \dbE_t\Big[|X^{t, y, Z, u}_T - \xi|^2\Big],\q (t,y) \in [0, T] \times \dbR^{d'}.
\eeaa
Similar to \reff{Wyreg} and by the uniform boundedness in Assumption \ref{assum-ocont}, one can choose a version of $\tilde W$ such that 
\bea
\label{tildeWreg}
|\tilde W(t, y_1)-\tilde W(t, y_2)|\le  C[1+|y_1|+|y_2|]|y_1-y_2|, \q\dbP_0\mbox{-a.s.}
\eea
 Then by standard arguments one can easily show that
\bea
\label{tildeW}
\tilde W(t, \eta) = \essinf_{(Z,u)\in  \hL^2(\hF, \hR^{d'\times d}) \times \sU} \dbE_t\Big[|X^{t, \eta, Z, u}_T - \xi|^2\Big],\q \forall \eta \in \dbL^2(\cF_t, \dbR^{d'}).
\eea
Next, following the arguments in Theorem \ref{thm-duality}, one can prove the following duality results:

(D-i) If  $u\in \sU$, and $(Y^u,Z^u)$ is the solution to BSDE (\ref{Y1}), then $\tilde W(t, Y^u_t) = 0$, $\hP_0$-a.s.;

(D-ii)   If $\eta \in \hL^0(\cF_t, \dbR^{d'})$ is such that $\tilde W(t, \eta)=0$, $\hP_0$-a.s., then $\eta \in \hL^2(\cF_t, \hR^{d'})$. 
Furthermore, for any $\e>0$,  there exists $u^\e\in \sU$, such that
\bea
\label{tildeDuality}
 |Y^{u^\e}_t-\eta|\le C\e, \q \mbox{and} \q \lim_{\e\to 0} \Phi(t, Y^{u^\e}_t) = \Phi(t, \eta),~\dbP_0\mbox{-a.s.}
\eea

We now construct the family of maximizers $\{\ol Y_t\}$.  For each fixed $t\in[0,T]$, denote 
$\tilde \sN_t := \{(\o, y)\in \O\times \hR^{d'}: \tilde W(t,\o, y) =0\}$ and $\tilde \sN_t(\o) := \{y\in \hR^{d'}: (\o, y)\in \tilde\sN_t \}$.
Then $\tilde \sN_t$ is $\cF_t\times \cB(\hR^{d'})$-measurable, and for $\hP_0$-a.e. $\o\in\O$, $\tilde \sN_t(\o)$ is closed and bounded, whence compact.  Now define $\ol \Phi_t(\o) := \sup_{y\in \tilde \sN_t(\o)} \Phi(t,\o, y)$, and denote
\beaa
\sM_t := \{(\o, y)\in \tilde \sN_t: \ol \Phi_t(\o) = \Phi(t,\o, y)\},\qq \sM_t(\o):= \sM_t \cap \tilde \sN_t(\o).
\eeaa
Then it is easy to see that $\ol \Phi_t$ is $\cF_t$-measurable and $\sM_t$ is $\cF_t \times \cB(\hR^{d'})$-measurable. Moreover,
the continuity of $\Phi$ in $y$ implies that $\sM_t(\o)$ is nonempty and compact, for $\hP_0$-a.e. $\o\in\O$. Now let $\ol Y_t(\o)$ be the (unique) maximum point of $\sM_t(\o)$ under the following order on $\hR^{d'}$:
\beaa
y < y' \q\Longleftrightarrow\q \mbox{for some} ~i=1,\cds, d',  \q y_j = y'_j, ~j=1,\cds, i-1, ~\mbox{and}~ y_i < y'_i.
\eeaa
Then clearly $\ol Y_t$ is $\cF_t$-measurable, and $\ol Y_t(\o) \in \cM_t(\o)$.

We now verify that $\ol Y$ satisfies all the requirements in Definition \ref{defn-consistent}. First, it is clear that $\tilde\sN_T(\o) = \{\xi(\o)\}$, and thus $\ol Y_T(\o) = \xi(\o)$. We next show that 
\bea
\label{PhiolY}
\ol Y_t \in \dbL^2(\cF_t, \dbR^{d'})\q\mbox{and} \q \Phi(t, \ol Y_t) = \esssup_{u\in \sU} \Phi(t, Y^u_t),
\eea
Indeed, for any $u\in \sU$, by the duality result (D-i) above we  have $\tilde W(t, Y^u_t) = 0$. That is, $Y^u_t(\o) \in \tilde \sN_t(\o)$, and thus $\Phi(t, \o, Y^u_t(\o)) \le \ol \Phi_t(\o) = \Phi(t,\o, \ol Y_t(\o)$, for $\dbP_0$-a.e. $\o\in \O$. Conversely, since $\ol Y_t(\o) \in \sM_t(\o) \subset \tilde\sN_t(\o)$, we see that $\tilde W(t, \o, \ol Y_t(\o))=0$ for $\dbP_0$-a.e. $\o$.  Then by the duality result (D-ii) and \reff{tildeDuality} we prove \reff{PhiolY} immediately.

It remains to verify the DPP \reff{olYDPP}. Note that for any $u^0\in \sU$,  \reff{PhiolY} implies that $\Phi(t_2, Y^{u^0}_{t_2}) \le \Phi(t_2, \overline Y_{t_2})$, $\dbP_0$-a.s. Then, it follows from the comparison principle 
\reff{Phicomparison} that
\beaa
\esssup_{u\in\sU} \Phi(t_1, \sY^u_{t_1}(t_2, Y^{u^0}_{t_2}) \le \esssup_{u\in\sU} \Phi(t_1, \sY^u_{t_1}(t_2, \overline Y_{t_2})), \qq\dbP_0\mbox{-a.s.}
\eeaa
Note that by definition $\sY^u_{t_1}(t_2, Y^{u^0}_{t_2}) = Y^{u\otimes_{t_2} u^0}_{t_1}$, then clearly
\beaa
\Phi(t_1, \overline Y_{t_1})= \esssup_{u\in\sU} \Phi(t_1, Y^u_{t_1}) \le \esssup_{u\in\sU} \Phi(t_1, \sY^u_{t_1}(t_2, \overline Y_{t_2})),~\dbP_0\mbox{-a.s.}
\eeaa
On the other hand, again by \reff{PhiolY}, there exist $\{u^\e\}_{\e>0}\subseteq \sU$ such that $ |Y^{u^\e}_{t_2}-\ol Y_{t_2}|\le C\e$, $\dbP_0$-a.s.
Then for any $u\in \sU$, by the stability of BSDE and the continuity of $\Phi$ in $y$,
\beaa
\Phi(t_1, \sY^u_{t_1}(t_2, \overline Y_{t_2})) = \lim_{\e\to 0} \Phi(t_1, \sY^u_{t_1}(t_2, Y^{u^\e}_{t_2})) = \lim_{\e\to 0} \Phi(t_1, Y^{u\otimes_{t_2}u^\e}_{t_1}) \le \Phi(t_1, \overline Y_{t_1}).
\eeaa
Since $u\in \sU$ is arbitrary, we obtain \reff{olYDPP}, completing the proof.     
\qed

\subsection{The linear case}

While Theorem \ref{thm-comparison} gives a guiding principle for finding the time consistent dynamic utility function, it would be
extremely desirable to see if a function satisfying the comparison principle \reff{Phicomparison} does exist. In this subsection we shall construct an explicit example, in the case 
when both  BSDE (\ref{Y1}) and function $\f$ are linear. Our construction follows the dimension reduction technique in \cite{KZ}.   
\begin{thm}
\label{thm-linear}
Let Assumption \ref{assum-f} hold and assume that the coefficients $f$ and $\f$ are of the following linear form: 
\bea
\label{flinear}
\left.\ba{lll}
\dis f_i(t,\o, y, z, u) = \sum_{j=1}^{d'} [\a^{i,j}_t(\o) y_j + \b^{i,j}_t(\o) \cd z_j] + c_i(t,\o, u), \q i=1,\cds, d', \ms\\
\dis \f(y) = \sum_{i=1}^{d'} a_i y_i,
\ea\right.
\eea
Then there exists a random field $\Phi$ satisfying the  comparison principle \reff{Phicomparison}, 
which takes the following linear form:
\bea
\label{Philinear}
\Phi(t,\o,y) := \sum_{i=1}^{d'} A^i_t(\o) y_i,\q\mbox{with}~ A^i_0 = a_i,
\eea
\end{thm}
\proof  We first note that if $d'=1$, then the BSDE (\ref{Y1}) is 1-dimensional, thus the comparison theorem holds. Further
since $\f$ is linear, whence monotone, thus the problem is time consistent and the theorem becomes trivial. We shall thus concentrate on multi-dimensional cases. 
Note also that for $d'\ge 2$, following an inductional arguments as illustrated in \cite[Section 4.1]{KZ},
we need only prove the case $d'=2$. We shall split the proof (assuming $d'=2$) in three steps.

{\it Step 1.}  We begin by a heuristic argument which will lead us to the desired properties of the processes $A^1$ and $A^2$.
For convenience we shall assume that $A^1$ and $A^2$ take the form of  It\^o process:
\bea
\label{Ai}
A^i_t = a_i + \int_0^t b^i_s ds + \int_0^t \si^i_s dB_s,\q i=1,2,
\eea
For any $u\in\sU$ and the corresponding solution $(Y^u, Z^u)$, we define
\bea
\label{olYZu}
\hat Y^u_t := \Phi(t, \cd, Y^u_t):=\sum_{i=1}^2 A^i_t Y^{i,u}_t, \q \hat Z^u_t := \sum_{i=1}^2[A^i_t Z^{i,u}_t + \si^i_t Y^{i,u}_t], 
\q  t\in [0,T].
\eea
We hope to find a pair of processes  $(A^1, A^2)$ so that $(\hat Y^u, \hat Z^u)$  satisfy a one dimensional BSDE, so as to reduce the problem to the case $d'=1$. 

To this end, we first assume   $A^2_t \equiv a_2 \neq 0$, $0\le t\le T$, Then, an easy application of It\^o's formula and some direct computations lead us to
\bea
\label{dolYu}
 d \hat Y^u_t 
&=&\Big[A^1_t dY^{1,u}_t + Y^{1,u}_t dA^1_t  + \si^1_t Z^{1,u}_t dt + a_2 dY^{2,u}_t \Big]\\
&=& -\Big[A^1_t \sum_{j=1}^2[\a^{1,j}_t Y^{j,u}_t + \b^{1,j}_t Z^{j, u}_t] +A^1_t c_1(t, u_t) + a_2 \sum_{j=1}^2[\a^{2,j}_t Y^{j,u}_t + \b^{2,j}_t Z^{j, u}_t]  \nonumber\\
&&+a_2 c_2(t, u_t) - [b^1_t Y^{1,u}_t  +\si^1_t Z^{1,u}_t ] \Big]dt  +[A^1_t Z^{1,u}_t + \si^1_t Y^{1,u}_t+a_2Z^{2,u}_t ]dB_t. \nonumber\eea
Note that in this case $b^2=\si^2=0$, we see from (\ref{olYZu}) that   $A^1_t Z^{1,u}_t + \si^1_t Y^{1,u}_t+a_2Z^{2,u}_t =\hat Z^u_t $, and thus
\beaa
Y^{2,u}_t = a_2^{-1}[\hat Y^u_t -  A^1_t  Y^{1,u}_t],\q Z^{2,u}_t = a_2^{-1}[\hat Z^u_t - \si^1_t Y^{1,u}_t -  A^1_t  Z^{1,u}_t].
\eeaa
Plugging these into (\ref{dolYu}) and reorganizing terms  yields:
\bea
\label{dolYu1}
&&-d\hat Y^u_t+\hat Z^u_t dB_t\nonumber \\
&=&\Big[[A^1_t \a^{1,1}_t + a_2 \a^{2,1}_t - b^1_t] Y^{1,u}_t +  [A^1_t \b^{1,1}_t + a_2 \b^{2,1}_t-\si^1_t ] Z^{1,u}_t + A^1_t c_1(t, u_t) + a_2 c_2(t, u_t) \nonumber\\
&& + [A^1_t \a^{1,2}_t + a_2 \a^{2,2}_t]   a_2^{-1}[\hat Y^u_t -  A^1_t  Y^{1,u}_t]  +  [A^1_t \b^{1,2}_t + a_2 \b^{2,2}_t] a_2^{-1}[\hat Z^u_t - \si^1_t Y^{1,u}_t -  A^1_t  Z^{1,u}_t] \Big]dt  \nonumber\\
&=& \Big[ a_2^{-1} [A^1_t \a^{1,2}_t + a_2 \a^{2,2}_t]   \hat Y^u_t  + a_2^{-1} [A^1_t \b^{1,2}_t + a_2 \b^{2,2}_t] \hat Z^u_t+ A^1_t c_1(t, u_t) + a_2 c_2(t, u_t) \nonumber\\
&&+\Th_t Y^{1,u}_t +\G_t Z^{1,u}_t \Big] dt,
\eea
where 
\beaa
\Th_t&:=&A^1_t \a^{1,1}_t + a_2 \a^{2,1}_t - b^1_t ] - a_2^{-1}A^1_t [A^1_t \a^{1,2}_t + a_2 \a^{2,2}_t] - a_2^{-1} \si^1_t  [A^1_t \b^{1,2}_t + a_2 \b^{2,2}_t]; \nonumber\\
\G_t&:=&[A^1_t \b^{1,1}_t + a_2 \b^{2,1}_t-\si^1_t ]  - a_2^{-1} A^1_t [A^1_t \b^{1,2}_t + a_2 \b^{2,2}_t]. 
\eeaa
%
%

 Now setting $\Th_t\equiv \G_t\equiv 0$, we see that (\ref{dolYu1}) becomes a linear BSDE for $(\hat Y^u, \hat Z^u)$. 
But  this can be done by simply solving 
\beaa
b^1_t &:=& [A^1_t \a^{1,1}_t + a_2 \a^{2,1}_t ]  - a_2^{-1}A^1_t [A^1_t \a^{1,2}_t + a_2 \a^{2,2}_t]  - a_2^{-1} \si^1_t  [A^1_t \b^{1,2}_t + a_2 \b^{2,2}_t]; \\
\si^1_t &:=& [A^1_t \b^{1,1}_t + a_2 \b^{2,1}_t ]  - a_2^{-1} A^1_t [A^1_t \b^{1,2}_t + a_2 \b^{2,2}_t].
\eeaa
Note that the processes $b^1$ and $\si^1$ can be easily written as functions of the process $a_2^{-1}A^1$ by setting
 $b^1_t =   a_2 \hat b_1(t, \o, a_2^{-1} A^1_t)$ and $\si^1_t =  a_2 \hat\si_1(t, \o, a_2^{-1} A^1_t)$, where
\bea
\label{Ai=0}
\hat b_1(t,  x) &:=& |\b^{1,2}_t|^2 x^3 - \big[\a^{1,2} + \b^{1,2}[\b^{1,1}-\b^{2,2}] - \b^{1,2}\b^{22}\big] x^2\nonumber\\
&& + \big[\a^{1,1}_t - \a^{2,2}_t - \b^{2,2}[\b^{1,1}_t-\b^{2,2}_t]-\b^{1,2}_t\b^{2,1}_t\big] x + [\a^{2,1} - \b^{2,1}_t\b^{2,2}_t]; \\
\hat\si_1(t, x) &:=&  - \b^{1,2}_t |x|^2+ [\b^{1,1}_t-\b^{2,2}_t] x + \b^{2,1}_t. \nonumber
\eea
Plugging this into \reff{Ai}, we obtain an SDE for $A^1_t$:
\bea
\label{Ricatti}
A^1_t\slash a_2 = a_1 \slash a_2+ \int_0^t  \hat b_1(s, a_2^{-1} A^1_s) ds + \int_0^t  \hat \si_1(s, a_2^{-1} A^1_s) dB_s, \q t\ge 0.
\eea
We should note  that since the coefficients $\hat\si$ has quadratic growth in $A^1_t$ and $\hat b$ has triple growth in $A^1_t$, the
SDE \reff{Ricatti} is a {\it Ricatti equation} in general sense and has only local solutions. However, if (\ref{Ricatti}) is solvable, which we shall argue rigorously in the next step, 
then we will see that the $\Phi(t, \cd)$ defined by \reff{Philinear} satisfies the comparison principle \reff{Phicomparison}.

\ms
 
{\it Step 2.}  We now substantiate the idea in Step 1 rigorously. If $a_1 = a_2=0$, then clearly $V_0(\xi) = 0$ and there is nothing to prove. From now on we  assume without loss of generality that $|a_1|\le |a_2|$ and $a_2\neq 0$.  Denote $\t_0 := 0$. Recall \reff{Ricatti} and consider the following SDE:
\bea
\label{olAtruncate}
\hat A^1_t = a_1\slash a_2 + \int_0^t \hat b_1\big(s,  [-2] \vee \hat A^1_s \wedge 2\big) ds + \int_0^t \hat \si_1\big(s, [-2] \vee \hat A^1_s \wedge 2\big) dB_s, \q t\in [0,T].
\eea
Clearly $\hat A^1$ has global solution. Define
$\t_1 := \inf\{t\ge 0: |\hat A^1_t| \ge 2\} \wedge T$.
Then
\bea
\label{olA1}
\hat A^1_t = a_1\slash a_2 + \int_0^t \hat b_1\big(s,   \hat A^1_s\big) ds + \int_0^t \hat \si_1\big(s,  \hat A^1_s\big) dB_s,\q\t_0\le t\le \t_1.
\eea
We now set 
$A^1_t := a_2 \hat A^1_t$ and $A^2_t := a_2$, for  $\t_0 \le t\le \t_1$. 
Then, noting that $|\hat A^1_{\t_1}|=2$ (or $|(\hat A^1_{\t_1})^{-1}|={1\over 2}$) when $\t_1 < T$ and reversing the roles of $A^1$ and $A^2$ as in Step 1 we can then  obtain coefficients $\hat b_2$, $\hat \si_2$ completely symmetric as 
those in \reff{Ai=0}, and an SDE  on $[\t_1, T]$:
\beaa
\hat A^2_t = (\hat A^1_{\t_1})^{-1} + \int_{\t_1}^t \hat b_2\big(s,  [-2] \vee \hat A^2_s \wedge 2 \big) ds + \int_0^t \hat \si_2\big(s, [-2] \vee \hat A^2_s \wedge 2\big) dB_s.
\eeaa
Similarly $\hat A^2$ has global solution, and that
\bea
\label{olA2}
\hat A^2_t =(\hat A^1_{\t_1})^{-1} + \int_{\t_1}^t \hat b_2\big(s, \hat A^2_s\big) ds + \int_0^t \hat \si_2\big(s,  \hat A^2_s\big) dB_s,\q \t_1\le t\le \t_2,
\eea
where 
$\t_2 := \inf\{t\ge \t_1: |\hat A^2_t\slash A^1_{\t_1}| \ge 2\} \wedge T$.
We then define 
$A^1_t := A^1_{\t_1}$, and $A^2_t :=   A^1_{\t_1} \hat A^2_t$, for $\t_1 \le t\le \t_2$. 
Note that since
$A^1_{\t_1}  \hat A^2_{\t_1} = A^1_{\t_1}  (\hat A^1_{\t_1})^{-1} = a_2 = A^2_{\t_1}$, both $A^1$ and $A^2$ are continuous at 
$\t_1$.

Now repeating the arguments, we may define, for $n\ge 1$, processes $\{\hat A^n\}$ and stopping times $0=\t_0\le \t_1\le\t_n
\cds$,  such that
\beaa
&&\hat A^{2n}_t = (\hat A^{2n-1}_{\t_{2n-1}})^{-1} + \int_{\t_{2n-1}}^t \hat b_2\big(s,  \hat A^{2n}_s\big) ds + \int_{\t_{2n-1}}^t \hat \si_2\big(s, \hat A^{2n}_s\big) dB_s,\q \t_{2n-1}\le t\le \t_{2n};\\
&&\hat A^{2n+1}_t = (\hat A^{2n}_{\t_{2n}})^{-1} + \int_{\t_{2n}}^t \hat b_1\big(s,   \hat A^{2n+1}_s\big) ds + \int_{\t_{2n}}^t \hat \si_1\big(s,  \hat A^{2n+1}_s\big) dB_s,\q \t_{2n}\le t\le \t_{2n+1}.
\eeaa
Furthermore, for all $n\ge 1$, it holds that 
$|\hat A^{n}_t| < 2, ~\t_{n-1}\le t < \t_{n}$, and $|\hat A^{n}_{\t_{n}}| = 2$ on $\{\t_{n} <T\}$.
The rest of the argument will be based on the following fact, which will be validated in the next step:
\bea
\label{claim}
\dbP_0\Big(\bigcup_{n\ge 1}\{\t_n  = T\}\Big) = 1.
\eea

Assuming (\ref{claim}), we can now define continuous processes $A^1, A^2$ on $[0, T]$:
\bea
\label{A12}
\left.\ba{lll}
 A^1_t := A^1_{\t_{2n-1}},\q A^2_t := A^1_{\t_{2n-1}} \hat A^{2n}_t,,\q \t_{2n-1} \le t\le \t_{2n};\\
A^1_t := A^2_{\t_{2n}} \hat A^{2n+1}_t,\q A^2_t := A^2_{\t_{2n}},\q  \t_{2n} \le t\le \t_{2n+1}.
\ea\right.
\eea
Now define $\Phi$ by \reff{Philinear} and $(\hat Y^u, \hat Z^u)$ by \reff{olYZu}. We can rewrite \reff{dolYu} as
\beaa
d \hat Y^u_t  &=& -\Big[ \hat\a_t  \hat Y^u_t  + \hat\b_t \hat Z^u_t+ \sum_{i=1}^2A^i_t c_i(t, u_t) \Big]dt + \hat Z^u_t dB_t, \q 0\le t\le T,
\eeaa
where 
\bea
\label{alphabeta}
\left.\ba{lll}
\dis \hat\a_t&=&\left\{\ba{lll}
\a^{1,2}_t  \hat A^{2n+1}_t +  \a^{2,2}_t, \qq~&\mbox{on} ~[\t_{2n}, \t_{2n+1}]\\
\a^{2,1}_t  \hat A^{2n}_t+  \a^{1,1}_t, & \mbox{on}~ [\t_{2n-1}, \t_{2n}];
\ea\right.
\ms \\
\dis \hat\b_t&=&\left\{\ba{lll}
\b^{1,2}_t \hat A^{2n+1}_t +  \b^{2,2}_t, \qq~&\mbox{on} ~[\t_{2n}, \t_{2n+1}]\\
\b^{2,1}_t \hat A^{2n}_t  +  \b^{1,1}_t, & \mbox{on}~ [\t_{2n-1}, \t_{2n}].
\ea\right.
\ea\right. 
\eea
Note that $| \hat A^{2n+1}_t |\le 2$ on   $\t_{2n} \le t\le \t_{2n+1}$ and $|\hat A^{2n}_t |\le 2$ on   $\t_{2n-1} \le t\le \t_{2n}$, 
both $\hat\a, \hat\b$ are bounded. Now denoting $\hat Y^u_t(\xi)$ to emphasize the dependence on the terminal condition $\xi$, it follows from  the definition \reff{olYZu} and the comparison of BSDEs that 
\beaa
\Phi(T,\xi) \le \Phi(T,\tilde \xi)  &\Longrightarrow& \hat Y^u_t(\xi) \le \hat Y^u_t(\tilde \xi), ~\forall u\in\sU\\
&\Longrightarrow& \esssup_{u\in\sU} \Phi(t, \sY^u_{t}(T, \xi)) \le \esssup_{u\in\sU} \Phi(t, \sY^u_t(T, \tilde\xi)),~\hP_0\mbox{-a.s.}
\eeaa
The same argument can be used to treat any subinterval $[t_1, t_2]$, proving  \reff{Phicomparison}.

\ms
{\it Step 3.} It remains to prove \reff{claim}.  Fix some $\d>0$. Note that $|a_1\slash a_2|\le 1$. By \reff{olAtruncate} and standard estimates for SDEs we can easily check that
$\hE\big[\sup_{0\le t\le T} |\hat A^1_t|^2\big] \le C$.
Thus
\beaa
&&\dbP_0(\t_1 < T\wedge \d)  \le \dbP_0\Big(\sup_{0\le t\le \d} |\hat A^1_t| \ge 2\Big)\\
&\le& \dbP_0\Big(\sup_{0\le t\le \d} |\hat A^1_t - \hat A^1_0| \ge 1\Big) \le \dbE\Big[\sup_{0\le t\le \d}|\hat A^1_t - \hat A^1_0|^2\Big]\\
&\le& C\dbE\Big[\int_0^\d|\hat b_1\big(s,  [-2] \vee \hat A^1_s \wedge 2\big)|^2 ds + \int_0^\d |\hat \si_1\big(s, [-2] \vee \hat A^1_s \wedge 2\big)|^2 ds\Big]\le C\d.
\eeaa
Now setting $\d := {1\over 2C}$, so that 
\bea
\label{t1est}
\hP_0(\t_1 < T, \t_1 \le \d) \le {1\over 2}.
\eea
Similarly, noting that $|\hat A^2_{\t_1}| = {1\over 2}$ and $|\hat A^2_{\t_2}| = 2$ on $\{\t_2 < T\}$, we have 
\bea
\label{t2est}
\hP_0\Big(\t_2 < T\wedge (\t_1+ \d)~\Big|\cF_{\t_1}\Big) \le {1\over 2}.
\eea
Repeating the arguments, for any $n$ one shows that
\bea
\label{tnest}
\dbP_0\Big(\t_{n+1} < T\wedge (\t_n + \d)~\Big|\cF_{\t_n}\Big) \le {1\over 2}.
\eea

We shall prove (\ref{claim}) by arguing that $\hP_0\big\{\big(\bigcup_{n\ge 1}\{\t_n=T\}\big)^c\big\}=\hP_0\big\{\bigcap_{n\ge1}\{
\t_n<T\}\big\}=0$. But since $\t_n$'s are increasing, this amounts to saying that $\lim_{n\to\infty} \hP_0\{\t_n<T\}=0$. Now for
the given $\d$, we can assume that $m\d < T \le (m+1)\d$, for some $m\in\hN$. We claim the following much stronger
result,  which obviously  implies  \reff{claim}: for any $n\ge 1$, 
\bea
\label{induction}
\dbP_0(\t_n < T)  \le {(2n)^m\over 2^n}, \qq \mbox{whenever}  ~~m\d <T\le (m+1)\d.
\eea
We shall prove (\ref{induction}) by induction on $m$.  First, if $m=0$, namely $0< T \le \d$, then
\beaa
\dbP_0(\t_n < T) &=& \dbP_0(\t_n< T, \t_1 \le \d) = \dbP_0(\t_1 < T, \t_1 \le \d) \dbP_0\Big(\t_n < T\Big|\cF_{\t_1}, \t_1 < T\Big)\\
& \le& {1\over 2}\dbP_0\Big(\t_n < T\Big|\cF_{\t_1}, \t_1 < T\Big).
\eeaa
thanks to \reff{t1est}. By \reff{tnest}, for $k<n$ we have
\beaa
\dbP_0\Big(\t_n < T\Big|\cF_{\t_{k-1}}, \t_{k-1} < T\Big) \le {1\over 2}\dbP_0\Big(\t_n < T\Big|\cF_{\t_k}, \t_k < T\Big).
\eeaa
Then  by induction we  see that
\beaa
\dbP_0(\t_n < T) & \le& {1\over 2^{n-1}}\dbP_0\Big(\t_n < T\Big|\cF_{\t_{n-1}}, \t_{n-1} < T\Big)\le {1\over 2^n},
\eeaa
proving \reff{induction} for $m=0$. 

Assume \reff{induction} holds for $m-1$ and we shall prove it for $m$.  By \reff{t1est} we have
\beaa
\dbP_0(\t_n < T) &=& \dbP_0(\t_n< T, \t_1 \le \d) + \dbP_0(\t_n < T, \t_1 >\d) \\
&\le&  \dbP_0(\t_1 < T, \t_1 \le \d) \dbP_0\Big(\t_n < T\Big|\cF_{\t_1}, \t_1 < T\Big) + \dbP_0(\t_n < T, \t_n -\t_1 < T-\d)\\
& \le& {1\over 2}\dbP_0\Big(\t_n < T\Big|\cF_{\t_1}, \t_1 < T\Big) +  \dbP_0(\t_n < T, \t_n -\t_1 < T-\d).
\eeaa
Note that $(m-1)\d < T-\d \le m\d$, then the inductional hypothesis implies that 
\beaa
\dbP_0(\t_n < T, \t_n -\t_1 < T-\d) \le {(2n-2)^{m-1}\over 2^{n-1}},
\eeaa
and thus
\beaa
\dbP_0(\t_n < T) & \le& {1\over 2}\dbP_0\Big(\t_n < T\Big|\cF_{\t_1}, \t_1 < T\Big) +  {(2n-2)^{m-1}\over 2^{n-1}}.
\eeaa
By \reff{tnest} , for $k<n$ we have
\beaa
\dbP_0\Big(\t_n < T\Big|\cF_{\t_{k-1}}, \t_{k-1} < T\Big) & \le& {1\over 2}\dbP_0\Big(\t_n < T\Big|\cF_{\t_k}, \t_k < T\Big) +  {(2n-2k)^{m-1}\over 2^{n-k}}.
\eeaa
Then by induction we have
\beaa
\dbP_0(\t_n < T) \le {1\over 2^n} + \sum_{k=1}^{n-1} {(2k)^{m-1}\over 2^{n-1}}  = {1+ 2\sum_{k=1}^{n-1} (2k)^{m-1}\over 2^n}.
\eeaa
It is straightforward to check that
$1+ 2\sum_{k=1}^{n-1} (2k)^{m-1} \le (2n)^m$, proving \reff{claim}, whence the theorem. 
\qed

\section{The Master Equation Approach}
\label{sect-master}
\setcounter{equation}{0}

In this section we deviate from the dynamic utility $\Phi$ and attack the value function $V_0(\xi)$ from a different direction. 
We begin by noticing that, unlike
the forward stochastic control problem where the value function depends on the ``initial data", in our problem
the value $V_0(\xi)$ should be considered as the function of terminal data $(T, \xi)$. Our main idea is to let
$(T, \xi)$  become ``variables", and study the behavior of the value function.   For notational simplicity, in this section we denote $\dbL^2(\cF_t) := \dbL^2(\cF_t, \dbR^{d'})$.

To be more precise, let us consider the following set
\bea
\label{Th}
\sA :=\Big\{(t, \eta): t\in [0, T], \eta  \in \dbL^2(\cF_t)\Big\} \subset [0, T] \times \dbL^2(\cF_T).
\eea
We should note that the pair $(t, \eta)\in\sA$ is ``progressively measurable" in nature, that is, for 
each $t$, $\eta$ has to be  $\cF_t$-adapted. 

We now introduce a dynamic ``value" function for our original problem. Let  $\Psi: \sA \to \hR$ be a real-valued
function on $\sA$ defined by 
\bea
\label{Psi}
\Psi(t,\eta) =  \sup_{u\in \sU} \f(\sY^u_0(t, \eta)),\q (t,\eta)\in \sA.
\eea 
Clearly, it holds that
\bea
\label{Psi0}
\Psi(0, y)= \f(y) \q\mbox{and}\q V_0(\xi) = \Psi(T,\xi).
\eea

Furthermore, we have the following easy consequences for the value function $\Psi$. Among other things, we 
show that a ``forward" dynamic programming principle actually holds without any extra conditions, even in such 
a time-inconsistent setting. 
\begin{lem}
\label{lem-PsiDPP}
Assume that Assumption \ref{assum-f} is in force. Then, 

(i) For each $t$, $\Psi(t,\cd): \dbL^2(\cF_t)\to \dbR$ is Lipschitz  continuous: 
\bea
\label{Lipschitz}
|\Psi(t, \eta_1)-\Psi(t, \eta_2)| \le C\|\eta_1-\eta_2\|_{\dbL^2(\cF_t)}\q\mbox{for any}~\eta_1, \eta_2\in \dbL^2(\cF_t).
\eea

(ii) $\Psi$ satisfies the following ``forward dynamic programming principle": 
\bea
\label{PsiDPP}
\Psi(t_2, \eta) = \sup_{u\in \sU} \Psi(t_1, \sY^u_{t_1}(t_2, \eta)),\q \forall 0\le t_1 < t_2\le T, \eta\in \dbL^2(\cF_{t_2}).
\eea
\end{lem}
\proof (i) For any $\eta_1, \eta_2\in \dbL^2(\cF_t)$ and any $u\in \sU$, by standard BSDE arguments we have
\beaa
|\sY^u_0(t, \eta_1) - \sY^u_0(t, \eta_2)|^2 \le C\hE[|\eta_1-\eta_2|^2].
\eeaa
This immediately leads to  \reff{Lipschitz} since $u \in\sU$ is arbitrary.

(ii) Let $u\in \sU$ be given. By the uniqueness of the BSDE 
we should have
\beaa
\f\big(\sY^u_0(t_2, \eta)\big) = \f\Big(\sY^u_0\big(t_1, \sY^u_{t_1}(t_2, \eta)\big)\Big) \le \Psi\big(t_1,   \sY^u_{t_1}(t_2, \eta)\big).
\eeaa
 Taking supremum over $u$ we prove ``$\le$" part of \reff{PsiDPP}.  To see the opposite inequality, we fix an arbitrary $u\in \sU$.  For any $\e>0$, by the definition of $\Psi$, there exists $u^\e \in \sU$ such that
 \beaa
 \Psi\big(t_1,   \sY^u_{t_1}(t_2, \eta)\big) \le \f\Big(\sY^{u^\e}_0\big(t_1, \sY^u_{t_1}(t_2, \eta)\big)\Big)  + \e =  \f\Big(\sY^{u^\e\otimes_{t_1} u}_0(t_2, \eta)\Big)  + \e\le  \Psi(t_2, \eta) + \e.
 \eeaa
 Taking supremum over $u\in\sU$ on left side and sending $\e$ to zero in the right side, we obtain 
 the  ``$\ge$" part of \reff{PsiDPP} and completes the proof.
 \qed
 \begin{rem}
 \label{rem-forwardDPP}
 {\rm (i) Unlike the standard DPP in stochastic control literature, \reff{PsiDPP} is a forward DPP in the sense that the supremum in the right side acts on the smaller time $t_1$. This is due to the nature that our controlled dynamics is backward. This feature will also be crucial for deriving the master equation at below.
 
 (ii) In deterministic case, the $\Psi$ here coincides with the dynamic utility $\Phi$ constructed in \S\ref{sect-deterministic2}.
 \qed}
 \end{rem}
 
With the essentially ``free" dynamic programming equation (\ref{PsiDPP}), it is natural to envision an HJB-type equation
for the value function $\Psi$. We note that there are two fundamental differences between the current situation and the 
traditional ones: (i) since the DPP is ``forward", the HJB equation should also be a temporally forward
PDE; and (ii) since the spatial variable in the value function is now a random variable in an $\dbL^2$ space which is infinite
dimensional, the PDE is quite different from the traditional HJB equation (even those infinite dimensional ones(!)), due to its adaptedness requirement on the variable $\eta$. We therefore
call it {\it master equation}, which seems to fit the situation better than an ``HJB equation".

We now try to validate the idea. To begin with, we shall 
introduce appropriate notion of derivatives.  First, for each $t\in[0,T]$, viewing $\hL^2(\cF_t)$ as a Hilbert space and  denote by $\la\cd,\cd\ra$ its inner product, we can 
define the spatial derivative as the standard Fr\'echet derivative: for any $\eta, \tilde \eta\in \dbL^2(\cF_t)$,
\bea
\label{Deta0}
\la D_\eta \Psi(t, \eta), \tilde\eta\ra := \lim_{\e\to 0} {\Psi(t, \eta+\e\tilde\eta) - \Psi(t,\eta)\over \e},
 \eea
whenever the limit exists. We remark that, when $D_\eta \Psi(t, \eta)$ exists, it can (and will) be identified as a random 
variable in $\dbL^2(\cF_t)$, thanks to the Riesz Representation Theorem.
 
 The temporal derivative, however, is much more involved. We first note that the dynamic programming principle  \reff{PsiDPP} is ``forward", and more importantly, the value function is ``progressive measurable", it is conceivable that there might be some
 difference between two directional derivatives. As it turns out, if we use the right-derivative temporally as one often does, then the corresponding master equation will becomes obviously illposed. We shall provide a detailed analysis on this point in \S\ref{sect-right} below. 
 We will therefore use left-derivative.

A simple-minded, albeit natural, definition of the left-temporal derivative can be defined as follows:
 \bea
 \label{Dt1}
 \lim_{\d\to 0} {\Psi(t, \eta) - \Psi(t-\d,\eta)\over \d}.
 \eea
 However, bearing in mind the ``progressive measurability" of $\Psi$ (or the definition of the set $\sA$), we see that $\eta\in \dbL^2(\cF_t)$ is typically not $\cF_{t-\d}$-measurable, so $\Psi(t-\d,\eta)$ may not even be well-defined. One natural choice 
 to overcome this  issue  is to modify (\ref{Dt1}) to the following:
 \bea
 \label{Dt2}
 \lim_{\d\to 0} {\Psi(t, \eta) - \Psi(t-\d, \dbE^{\dbP_0}_{t-\d}[\eta])\over \d}.
 \eea
However, although this definition could actually be sufficient for our purpose in this paper, 
 it relies heavily on the underlying measure 
$\hP_0$, which would cause many unintended consequence when we encounter situations where various probability 
measures are involved, as we often see in applications. 

A  universal, ``measure-free", and potentially more applicable definition is the following ``pathwise" derivative:
 \bea
 \label{Dt}
 D^-_t \Psi(t,\eta) := \lim_{\d\to 0} {\Psi(t, \eta) - \Psi(t-\d,  \eta^t_{t-\d})\over \d}, ~\mbox{where}~ \eta^t_s(\o) := \eta(\o_{s\wedge \cd}), (s,\o)\in [0, t]\times \O.
 \eea
 provided the limit exists.  
 We remark that, $D^-_t \Psi(t,\eta)$ is a real number, if it exists.
 
 \ms
 Recall \S\ref{sect-path} for the notions in  pathwise analysis. We define
 
 \begin{defn} 
 \label{defn-C1}
 (i) $\Psi\in C^0(\sA)$ if $\Psi$ is continuous in $(t,\eta)$.

(ii) $\eta\in C^2(\cF_t)$  if the induced process $\eta^t\in C^{1,2}([0, t]\times \O)$. In this case, we denote
\bea
\label{paeta}
\pa_t \eta := \pa_t \eta^t_t,\q \pa_\o \eta := \pa_\o \eta^t_t,\q \pa^2_{\o\o} \eta := \pa^2_{\o\o} \eta^t_t.
\eea
Moreover, denote $C^2_b(\cF_t):= \{\eta \in C^2(\cF_t):  \eta, \pa_t \eta^t, \pa_\o \eta^t, \pa^2_{\o\o}\eta^t ~\mbox{are bounded}~\}$.

 (iii) $\Psi \in C^1(\sA)$ if $\Psi \in C^0(\sA)$, $D_\eta \Psi$ exists and is in $C^0(\sA)$, and  $D^-_t\Psi(t, \eta)$ exists for all $(t,\eta)\in \sA_0$, where
 \bea
 \label{sA0}
 \sA_0 := \{(t, \eta): 0\le t\le T, \eta \in C^2_b(\cF_t)\} \subset \sA.
 \eea
 \end{defn}
 
 \no We remark that, for $\Psi\in C^0(\sA)$, it is uniquely determined by its values in $\sA_0$.
 
 The main result of this section is the following theorem.
 \begin{thm}
 \label{thm-Master}
 Let Assumption \ref{assum-f} hold and $f(t,\o, 0,0, u)$ be bounded. Assume the $\Psi$ defined by \reff{Psi} is in  $C^1(\sA)$. Then, $\Psi$ satisfies the following master equation on $\sA$:
 \bea
 \label{Master}
 \left\{\ba{lll}
\dis D^-_t\Psi(t,\eta) =  \la D_\eta \Psi(t, \eta), \pa_t \eta + {1\over 2}\tr( \pa^2_{\o\o}\eta)\ra\\
\dis\qq +  \sup_{u\in \dbL^0(\cF_t, U)}
\la D_\eta \Psi(t, \eta), f(t, \eta, \pa_\o\eta, u)\ra,\q  (t,\eta)\in \sA_0;\\
 \dis \Psi(0, y) = \f(y),\q y\in \dbR^{d'}.
 \ea\right.
 \eea
 \end{thm}
 \proof Fix $0<\d<t$. We first apply the {\it functional It\^o formula} (\ref{FIto}) to get
 \beaa
 \eta^t_s = \eta - \int_s^t [\pa_t \eta^t_r + {1\over 2} \tr(\pa^2_{\o\o}\eta^t_r)] dr - \int_s^t \pa_\o \eta_r \cd dB_r, \q t-\d\le s\le t,~\dbP_0\mbox{-a.s.}
 \eeaa
For any $u\in \sU$, let $(\sY^u, \sZ^u):=(\sY^u(t, \eta), \sZ^u(t, \eta))$ be the solution to BSDE \reff{cY}.  Denote
 \beaa 
 \D Y^u_s :=  \sY^u_s - \eta^t_s,\q \D Z^u_s := \sZ^u_s - \pa_\o \eta^t_s, \q   t-\d\le s\le t. 
 \eeaa
 Then
 \bea
\label{DYu1}
\D Y^u_s = \int_s^t \Big[f(r, \sY^u_r, \sZ^u_r, u_r)  -  [\pa_t \eta^t_r + {1\over 2} \tr(\pa^2_{\o\o} \eta^t_r)] \Big]dr + \int_s^t \D Z^u_r dB_r,
  t-\d\le s\le t.
\eea
By standard BSDE estimates we see that
 \beaa
 \dbE\Big[\sup_{t-\d\le s\le t} |\D Y^u_s|^2 + \int_{t-\d}^t |\D Z^u_s|^2 ds \Big] \le C\d^2.
 \eeaa

We can now apply the forward dynamic programming principle \reff{PsiDPP} to get 
 \beaa
 \Psi(t,\eta) - \Psi(t-\d, \eta^t_{t-\d}) &=& \sup_{u\in \sU}\Big[\Psi\big(t-\d, \sY^u_{t-\d}\big) - \Psi(t-\d, \eta^t_{t-\d})\Big] \\
 &=& \sup_{u\in \sU} \int_0^1 \Big\la D_\eta \Psi\big(t-\d, \eta^t_{t-\d} + \th \D Y^u_{t-\d}\big), ~\D Y^u_{t-\d}\Big\ra d\th.
 \eeaa
To identify the right hand side above, we first deduce from (\ref{DYu1}) that
\beaa
I^u_\d &:=&  \D Y^u_{t-\d} - \int_{t-\d}^t \dbE_{t-\d}\Big[f(s, \eta^t_s, \pa_\o\eta^t_s, u_s)  -  [\pa_t \eta^t_s + {1\over 2} \tr(\pa^2_{\o\o} \eta^t_s)] \Big]ds\\
&=& \int_{t-\d}^t \dbE_{t-\d}\Big[ f(s, \sY^u_s, \sZ^u_s, u_s)  -  f(s, \eta^t_s, \pa_\o\eta^t_s, u_s) \Big]ds.
\eeaa
Then, it is not hard to check, using Assumption \ref{assum-f}, that
\beaa
\hE[|I^u_\d|^2] &\le& C\d\dbE\Big[\int_{t-\d}^t [|\D Y^u_s|^2+|\D Z^u_s|^2] ds\Big] \le C\d^3.
\eeaa
Consequently, as $\d\to 0$, we have
 \beaa
&& \Psi(t,\eta) - \Psi(t-\d, \eta^t_{t-\d})\\
&=&  \sup_{u\in \sU} \Big\la \int_0^1D_\eta \Psi\big(t-\d, \eta^t_{t-\d} +\th \D Y^u_{t-\d}\big)d\th, \\
&&\qq \int_{t-\d}^t \dbE_{t-\d}\Big[f(s, \eta^t_s, \pa_\o\eta^t_s, u_s)  -  [\pa_t \eta^t_s + {1\over 2} \tr(\pa^2_{\o\o} \eta^t_s)] \Big]ds+I^u_\d\Big\ra\\
 &=&  \sup_{u\in \sU} \Big\la D_\eta \Psi\big(t-\d, \eta^t_{t-\d} \big), \int_{t-\d}^t \dbE_{t-\d}\Big[f(s, \eta^t_s, \pa_\o\eta^t_s, u_s)  -  [\pa_t \eta^t_s + {1\over 2} \tr(\pa^2_{\o\o} \eta^t_s)] \Big]ds\Big\ra  + o(\d)\\
  &=&  \sup_{u\in\sU} \Big\la D_\eta \Psi\big(t-\d, \eta^t_{t-\d} \big), ~\int_{t-\d}^t \Big[f(s, \eta^t_s, \pa_\o\eta^t_s, u_s)  -  [\pa_t \eta^t_s + {1\over 2}\tr( \pa^2_{\o\o} \eta^t_s)] \Big]ds\Big\ra   + o(\d)\\
 &=& \sup_{u\in\sU}  \Big\la D_\eta \Psi\big(t, \eta \big), ~\int_{t-\d}^t \Big[f(t, \eta, \pa_\o\eta, u_s)  -  [\pa_t \eta + {1\over 2} \tr(\pa^2_{\o\o} \eta)] \Big]ds\Big\ra  + o(\d)\\
 &=&\d \sup_{u\in \dbL^0(\cF_t, U)}\Big\la D_\eta \Psi\big(t, \eta \big), ~f(t, \eta, \pa_\o\eta, u)  -  [\pa_t \eta + {1\over 2}\tr( \pa^2_{\o\o} \eta)]\Big\ra  + o(\d).
 \eeaa
 This implies \reff{Master} immediately.
 \qed
 
 \begin{rem}
 \label{rem-left}
 {\rm 
(i) From (\ref{Master}) we see that the master equation is a first order (forward) path-dependent PDE (although it involves 
 the second-order path-derivative of the state variable $\eta$). While this is 
 obviously the consequence of the forward DPP \reff{PsiDPP} and our required initial condition on $\Psi$, it is also 
 due to the fact that, for a forward problem, standing at $t$ and looking ``left", the problem is essentially "deterministic", 
hence the corresponding ``HJB equation should be first order. The {\it left-temporal path derivative} that we introduced
in (\ref{Dt}) is thus essential.


(ii) The main difficulty of this approach is the proper solution of the master equation (\ref{Master}). To the best of our knowledge, 
such an equation is completely new in the literature. Its wellposedness, in strong, weak, and viscosity sense, seem to be all 
open at this point. We hope to be able to address some of them in 
our future research.
 \qed}
 \end{rem}


\subsection{An ill-posed master equation} 
\label{sect-right}
We have emphasized at above the importance for using the  left-temporal derivative, given the fact that $\Psi$ satisfies a forward dynamic programming principle. In what follows we shall reinforce this point by explaining how a 
``traditional" right-temporal derivative could actually lead to a ill-posed master equation.

Let us suppose that we define a time-derivative in the following traditional way:
\bea
\label{D+t}
D^+_t \Psi(t,\eta) := \lim_{\d\downarrow 0} {\Psi(t+\d, \eta) - \Psi(t,\eta)\over \d},\q (t,\eta) \in \sA,
\eea
whenever the limit exists. Since by our definition of $\sA$, for each $\d>0$, $\eta\in \dbL^2(\cF_t) \subset \dbL^2(\cF_{t+\d})$, thus $\Psi(t+\d, \eta)$ 
is well-defined for all $(t,\eta)\in\sA$.

Now let us derive the equation for the $\Psi$ in \reff{Psi}  involving such a derivative. Again, by DPP \reff{PsiDPP} we have
\bea
\label{Psit+d}
\Psi(t+\d, \eta) - \Psi(t,\eta) &=& \sup_{u\in \sU} \Big[\Psi(t, \sY^u_t(t+\d, \eta)) - \Psi(t, \eta)\Big]\nonumber\\
& =& \sup_{u\in \sU} \int_0^1 \Big\la D_\eta \Psi(t, \eta + \th \cY^u_t), \cY^u_t\Big\ra d\th
\eea
where $\cY^u_s := \sY^u_s(t+\d, \eta) - \eta$, $t\le s\le t+\d$. Note that, if we denote $\cZ^u_s:= \sZ^u_s(t+\d, \eta)$, then
$(\cY^u, \cZ^u)$ satisfies the BSDE:
\beaa
\cY^u_s = \int_s^{t+\d}f(r, \eta+\cY^u_r, \cZ^u_r, u_r) dr - \int_s^{t+\d} \cZ^u_r dB_r, \q t\le s\le t+\d.
\eeaa
Then, the  standard BSDE estimates would tell us that,
\beaa
\dbE\Big[\sup_{t\le s\le t+\d} |\cY^u_s|^2 + \int_t^{t+\d}|\cZ^u_s|^2 ds\Big] \le C\d^2.
\eeaa
Again, let us denote
\beaa
I^u_\d&:=& \cY^u_t - \dbE_t\Big[\int_t^{t+\d} f(s, \eta, 0, u_s) ds\Big].
\eeaa
Then, assuming Assumption \ref{assum-f} we have
\beaa
|I^u_\d| = \Big|\dbE_t\Big[\int_t^{t+\d}[f(s, \eta+\cY^u_s, \cZ^u_s, u_s) - f(s, \eta, 0, u_s) ]ds\Big| \le C\dbE_t\Big[\int_t^{t+\d}[|\cY^u_s|+|\cZ^u_s|]ds\Big],
\eeaa
and consequently 
\beaa
\dbE[|I^u_\d|^2] \le C\d \dbE\Big[\int_t^{t+\d}[|\cY^u_s|^2+|\cZ^u_s|^2\Big]ds\le C \d^3.
\eeaa
Now \reff{Psit+d} will lead to that
\beaa
\Psi(t+\d, \eta) - \Psi(t,\eta)  & =& \sup_{u\in\sU} \Big\la D_\eta \Psi(t, \eta), ~\dbE_t\Big[\int_t^{t+\d} f(s, \eta, 0, u_s) ds\Big]\Big\ra + o(\d)\\
&=& \d \sup_{u\in \dbL^2(\cF_t, U)} \Big\la D_\eta \Psi(t, \eta), f(t, \eta, 0, u) \Big\ra  + o(\d).
\eeaa
In other words, we will arrive at the following first order PDE:
\bea
\label{Master+}
\left\{ \ba{lll} 
\dis D^+_t  \Psi(t,\eta) =\sup_{u\in \dbL^2(\cF_t,U)} \Big\la D_\eta \Psi(t, \eta), f(t, \eta, 0, u) \Big\ra, \q (t,\eta)\in \sA; \ms\\
\Psi(0, y) = \f(y).
\ea\right.
\eea
We remark that the equation \reff{Master+} is typically ill-posed. Indeed, \reff{Master+} involves only $f(\cd,\cd, 0, \cd)$, while the $\Psi$ defined in \reff{Psi} obviously depends on $f(\cd, \cd, z, \cd)$. So unless the function $f$ is independent of the variable $z$, 
there is essentially no hope that the  equation \reff{Master+} will have a unique solution, as the value functions of two completely
different optimization problems can satisfy the same master equation(!). 
We therefore conclude that $D^-_t \Psi$, not $D^+_t\Psi$, is the right choice of temporal derivative for the master equation.

\end{document}